\documentclass[10pt,leqno]{article}
\usepackage[utf8]{inputenc}

\date{}
\usepackage{mathrsfs}
\usepackage{hyperref} 
\newfont{\cyr}{wncyr10}

\usepackage[francais]{babel}
\usepackage{amssymb}
\date{}
\input{xypic}

\renewcommand{\leq}{\leqslant}
\renewcommand{\geq}{\geqslant}

\newcommand{\hotimes}{\widehat{\otimes}}
\renewcommand{\phi}{\varphi}

\renewcommand{\dim}[1]{\mbox{\rm dim}_{#1}\;}

\newcommand{\got}[1]{{\mathfrak #1}}

\renewcommand{\Bbb}{\mathbb}    
\renewcommand{\cal}{\mathscr}

\newcommand{\RR}{{\Bbb R}}
\newcommand{\KK}{{\Bbb K}}
\newcommand{\ZZ}{{\Bbb Z}}

\newcommand{\DD}{{\Bbb D}}

\newcommand{\Aff}{{\Bbb A}}
\newcommand{\NN}{{\Bbb N}}
\newcommand{\CC}{{\Bbb C}}

\newcommand{\QQ}{{\Bbb Q}}

\renewcommand{\epsilon}{\varepsilon}

\newcommand{\zero}{^{\mbox{\tiny o}}}

\newcommand{\spec}{\mbox{\rm Spec}\;}

\newcommand{\red}{\widetilde}
\newcommand{\gred}[1]{\widetilde{#1}_{\bullet}}

\newcommand{\hres}{{\cal H}}

\newcommand{\deux}[1]{\refstepcounter{subsection}\label{#1}\medskip\noindent {\bf (\thesubsection)}\hspace{.1cm}}
\newcommand{\trois}[1]{\refstepcounter{subsubsection}\label{#1}\medskip\noindent {\bf
    (\thesubsubsection)}\hspace{.1cm}}

\newcommand{\lr}{_{\leq r}}
\newcommand{\str}{_{<r}}
\newcommand{\hrr}{_{\prec r}}

\title{Toute forme modérément ramifiée d'un polydisque ouvert est triviale} 
\author{Antoine Ducros}

\begin{document}
\maketitle

{\small {\bf Abstract.} Let $k$ be a complete, non-Archimedean field and let $X$ be a $k$-analytic space. Assume that there exists a finite, tamely ramified extension $L$ of $k$ such that $X_L$ is isomorphic to an open polydisc over $L$ ; we prove that $X$ is itself isomorphic to an open polydisc over $k$. The proof consists in using the {\em graded} reduction (a notion which is due to Temkin) of the algebra of functions on $X$, together with some graded counterparts of classical commutative algebra results: Nakayama's lemma, going-up theorem, basic notions about \'etale algebras, etc.}

\section*{Introduction}

Le but de cet article est de démontrer le théorème suivant (th. \ref{pseudo}) : {\em soit $k$ un corps ultramétrique complet et soit $X$ un espace $k$-analytique. Supposons qu'il existe une extension finie, séparable et modérément ramifiée $L$ de $k$ telle que $X_{L}$ soit isomorphe à un $L$-polydisque ouvert de polyrayon $(r_1,\ldots,r_n)$ ; l'espace $X$ lui-même est alors isomorphe à un $k$-polydisque ouvert dont le polyrayon $(s_1,\ldots, s_n)$ est tel que $$|L^*|r_1^{\NN}r_2^{\NN}\ldots r_n ^{\NN}=|L^*|s_1^{\NN}s_2^{\NN}\ldots s_n ^{\NN}.$$ } 

\medskip
Notons un cas particulier important (corollaire \ref{corollnonramray1}) : si $X_L$ est isomorphe au $L$-polydisque {\em unité} ouvert de dimension $n$ et si $L$ est {\em non ramifiée} sur $k$, alors $X$ est isomorphe au $k$-polydisque unité ouvert de dimension $n$. Il suffit en effet pour le voir d'appliquer le théorème en supposant de surcroît que $|L^*|=|k^*|$ et que les $r_i$ sont tous égaux à $1$ ; les $s_i$ appartenant alors tous à $|k^*|$, le $k$-polydisque ouvert de polyrayon $(s_1,\ldots,s_n)$ est isomorphe au $k$-polydisque ouvert unité, d'où notre assertion. 

\subsection*{Quelques commentaires}

\begin{itemize}

\medskip
\item[1)] Lorsque nous parlons d'espace $k$-analytique, nous employons cette expression au sens de Berkovich (\cite{brk1}, \cite{brk2}). Toutefois l'essentiel de notre preuve est purement algébrique, et la nature précise des espaces en jeu n'y joue à peu près aucun rôle : tout se passe au niveau des anneaux de fonctions, qui sont insensibles au cadre de travail choisi (géométrie rigide, géométrie de Berkovich, géométrie de Huber, géométrie de Raynaud....). Les propriétés spécifiques à la théorie de Berkovich n'interviennent réellement qu'à deux propos : d'une part, parce que nous autorisons la valeur absolue de $k$ à être triviale (hypothèse exclue par les autres approches) ; et d'autre part lors d'un problème de descente, au cours de la preuve du théorème principal : il s'agit d'établir qu'un morphisme entre deux espaces sans bord\footnote{Cette assertion est vraie en général, {\em i.e.} avec bords éventuels, mais la preuve est alors plus délicate ; {\em cf.} par exemple \cite{cnrtmk}, th. 9.2. } qui devient un isomorphisme après extension des scalaires est déjà un isomorphisme sur le corps de base, ce que nous faisons en exploitant les bonnes propriétés topologiques des espaces de Berkovich, qui facilitent grandement les raisonnements locaux. 

\medskip
\item[2)] On peut démontrer le corollaire \ref{corollnonramray1} que nous évoquons ci-dessus (cas d'un polydisque unité et d'une extension non ramifiée) ou bien à partir de notre théorème principal, ou bien directement, en décalquant et simplifiant la preuve de ce dernier ({\em cf.} \ref{commentnonram}) ; expliquons succinctement les grandes étapes du raisonnement auquel conduirait cette démarche. 

\begin{itemize}
\medskip
\item[a)] Soit $X$ un espace $k$-analytique réduit, soit $\cal A$ l'algèbre des fonctions analytiques sur $X$ et soit $||.||_{\infty}$ la «norme spectrale» sur $X$ (qui peut prendre des valeurs infinies) ; soit ${\cal A}_{\leq 1}$ (resp. ${\cal A}_{<1}$)  le sous-anneau de $\cal A$ formé des fonctions $f$ telles que $||f||_{\infty}\leq 1$ (resp. $||f||_\infty<1$) et soit ${\cal A}_{\prec 1}$ le complété de ${\cal A}_{<1}$ pour la convergence uniforme sur les compacts. On note $\red {\cal A}$ le quotient ${\cal A}_{\leq 1}/{\cal A}_{\prec 1}$. Si $X$ est un polydisque unité ouvert de dimension $n$ alors $\red{\cal A}\simeq \red k [[\tau_1,\ldots, \tau_n]]$. 

\medskip
\item[b)] Si $L$ est une extension finie séparable et non ramifiée de $k$ alors pour toute extension complète $F$ de $k$ on a $\red{F\otimes_kL}\simeq \red F\otimes_{\red k}\red L$. 

\medskip
\item[c)] Soit $X$ un espace $k$-analytique et soit $L$ une extension finie séparable et non ramifiée de $k$ telle que $X_L$ soit isomorphe à un $L$-polydisque unité ouvert de dimension $n$ ; soit $\cal A$ (resp. $\cal B$) l'anneau des fonctions analytiques sur $X$ (resp. $X_L$) ; on montre à l'aide de b) que $\red {\cal B}$ (qui est isomorphe à $\red L [[\tau_1,\ldots, \tau_n]]$ en vertu de a) ) s'identifie naturellement à $\red{\cal A}\otimes_{\red k}\red L$. 

\medskip
\item[d)] À l'aide de c) et en exploitant le fait que $\red L$ est une extension séparable de $\red k$, on voit facilement que l'idéal maximal $\got m$ de $\red{\cal B}$ possède un système de $n$ générateurs appartenant à $\red{\cal A}$, que l'on peut relever en $n$ fonctions $f_1,\ldots, f_n$ dans $\cal A$. 

\medskip
\item[e)] Les $f_i$ induisent un morphisme $\phi$ de $X$ vers le polydisque unité ouvert de dimension $n$ sur $k$. On déduit du fait que les $\red {f_i}$ engendrent $\got m$ que $\phi_L$ est un isomorphisme ; par descente, $\phi$ est un isomorphisme. 

\end{itemize}

\medskip
\item[3)] Nous pouvons maintenant décrire comme suit la preuve de notre théorème principal : elle consiste {\em grosso modo} à reprendre les étapes a), b), c), d), et e) décrites ci-dessus, mais en rajoutant l'adjectif «gradué»~un peu partout. Plus précisément, nous utilisons une approche, introduite par Temkin dans \cite{tmk2}, qui consiste à étudier une algèbre normée $\cal A$ par le biais de son algèbre résiduelle {\em graduée}, qui est un anneau $\RR^*_+$-gradué, dont la partie de degré $1$ (attention, la graduation est multiplicative) est l'algèbre résiduelle classique. Le point clef sur lequel repose notre démarche est le suivant : si $k$ est un corps ultramétrique complet alors les extensions modérément ramifiées de $k$ se comportent, du point de vue de la théorie de la réduction graduée, exactement comme les extensions non ramifiées du point de vue de la réduction classique ; c'est très certainement, au moins dans le cas de la valuation discrète, une reformulation du fait que les extensions modérément ramifiées de $k$ sont précisément celles qui s'étendent en un log-schéma fini et log-étale sur $\spec k\zero$. 

\medskip
\item[4)] La motivation initiale de Temkin pour développer sa théorie était la suivante. La géométrie rigide ne considère que des polydisques fermés dont le polyrayon est constitué de réels appartenant à $\sqrt{|k^*|}$ ({\em i.e.} de réels qui sont de torsion modulo $|k^*|$), alors que la géométrie de Berkovich autorise tous les polyrayons. Or dès que le polyrayon d'un polydisque fermé $\DD$ d'algèbre des fonctions $\cal A$ contient un réel $r\in \RR^*_+-\sqrt{|k^*|}$, l'algèbre résiduelle classique $\red{\cal A}$ a un comportement inadéquat (par exemple, sa dimension de Krull est strictement inférieure à la dimension de $\DD$), et Temkin a réalisé que c'est son algèbre résiduelle {\em graduée} qui possède les propriétés requises. 

\medskip
Mais pour ce qui nous intéresse ici, à savoir les espaces qui sont des polydisques ouverts virtuels ({\em i.e.} après extension finie), nous ne pouvons utiliser la réduction classique {\em que} dans le cas des polydisques {\em unité} et des extensions non ramifiées ({\em cf.} le point 2) ci-dessus ) : dès que le polyrayon $(r_1,\ldots,r_n)$ de notre théorème principal contient au moins un réel $r_i$ qui n'appartient pas à $|k^*|$, ou dès que $|L^*|\neq |k^*|$, notre stratégie de preuve rend indispensable le recours au formalisme gradué, même si les $r_i$ appartiennent tous à $\sqrt {|k^*|}$, et bien que $|L^*|$ soit toujours contenu dans $\sqrt {|k^*|}$. 

\medskip
\item[5)] En ce qui concerne les polydisques {\em fermés}, l'analogue de notre théorème est faux, dès la dimension 1 : Tobias Schmidt explique dans \cite{tbs} comment construire un contre-exemple. Mais il y démontre tout de même que si $X$ est un espace $k$-analytique qui devient isomorphe au disque unité fermé sur une extension finie et {\em non ramifiée} de $k$, alors $X$ est déjà isomorphe au disque unité fermé sur $k$. 

\medskip
\item[6)] Cet article comprend trois parties. La première explique le formalisme de l'algèbre commutative graduée et énonce dans ce contexte un certain nombre d'avatars de résultats traditionnels, pratiquement sans la moindre preuve (certains raisonnements sont toutefois esquissés à propos de la théorie de Galois graduée) : les justifier avec précision reviendrait en effet à infliger au lecteur un long pensum sans grand intérêt, consistant en une retranscription {\em mutatis mutandis} des démonstrations classiques.

\medskip
Dans la seconde, après avoir rappelé en quoi consiste la théorie de la réduction graduée de Temkin, nous reformulons dans ce cadre la théorie de la ramification modérée. Les résultats que nous donnons concernant cette dernière sont bien connus, mais en raison de leur importance pour la suite, et du caractère inhabituel de la façon dont nous les présentons, nous avons choisi d'en redonner les preuves directement dans le langage de Temkin.

\medskip
 Quant à la troisième partie, elle est dévolue à la démonstration proprement dite du théorème principal, et se termine par quelques contre-exemples (\ref{contrex} {\em et sq.}) montrant que l'hypothèse de ramification modérée est indispensable à sa validité, puis par une question sur la validité éventuelle de l'analogue complexe de notre résultat. 
 
\end{itemize}

\section{Algèbre commutative graduée}\label{RAMI}

\medskip
Dans tout ce chapitre, on désigne par $\Gamma$ un groupe abélien divisible et sans torsion (autrement dit, un $\QQ$-espace vectoriel), que l'on note {\em multiplicativement} ; en particulier, l'élément neutre de $\Gamma$ est noté $1$. 

\medskip
{\em Remarque.} En pratique, $\Gamma$ sera égal à $\RR^*_+$. 

\subsection*{Algèbre commutative graduée} 

Les définitions, conventions et notations relatives à l'algèbre commutative graduée que nous utiliserons sont celles adoptées par Temkin dans \cite{tmk2} ; mais pour la commodité du lecteur, nous allons faire tous les rappels nécessaires. 

\deux{introalggrad} Un {\em anneau gradué} $A$ est un anneau commutatif muni d'une décomposition (comme groupe abélien) en somme directe $\bigoplus \limits_{\gamma \in \Gamma}  A_\gamma$ telle que $ A_\gamma. A_\delta\subset  A_{\gamma\delta}$ pour tout $(\gamma, \delta)$ ; on dira que $ A_\gamma$ est l'ensemble des {\em éléments homogènes de degré $\delta$.} Un idéal $ I$ de $ A$ est dit {\em homogène} s'il est engendré par des éléments homogènes ou, ce qui revient au même, s'il est égal à $\bigoplus ( I\cap  A_\gamma)$. Le quotient d'un anneau gradué par un idéal homogène hérite d'une graduation naturelle. Un morphisme d'anneaux gradués $ A\to  B$ est un morphisme d'anneau de $ A$ vers $ B$ qui envoie $ A_{\gamma}$ sur $ B_{\gamma}$ pour tout $\gamma$. 

\deux{trivgrad} Tout anneau commutatif peut être vu comme un anneau gradué par le biais de la {\em graduation triviale}, pour laquelle tout élément est homogène de degré 1, ce qui permet de voir l'algèbre commutative graduée comme une généralisation de l'algèbre commutative classique. Le but de ce qui suit est d'expliquer succinctement comment certaines définitions et théorèmes de cette dernière s'étendent à ce nouveau contexte ; nous ne donnerons pas le détail des démonstrations -- elles consistent en la reprise des preuves traditionnelles, agrémentées lorsqu'il convient de l'adjectif «homogène» ou «gradué» ; le lecteur pourra en trouver certaines dans \cite{tmk2} et \cite {cnrtmk}.

\trois{gradintegr} Soit $ A$ un anneau gradué. On dit qu'il est {\em intègre} s'il est non nul et si le produit de deux éléments {\em homogènes} non nuls de $ A$ est toujours non nul ; on dit que $ A$ est {\em réduit} si tout élément {\em homogène} et nilpotent de $ A$ est nul ; on dit que $ A$ est un {\em corps gradué} s'il est non nul et si tout élément {\em homogène} non nul de $ A$ est inversible. 

\trois{remcorpsgrad} {\em Remarque.} On prendra garde que l'anneau sous-jacent à un corps gradué n'est pas un corps en général ({\em cf. \ref{excorpsgrad} infra}) ; par contre, un anneau gradué est intègre (resp. réduit) au sens ci-dessus si et seulement si il est intègre (resp. réduit) comme anneau : la condition est en effet clairement suffisante, et on voit qu'elle est nécessaire en munissant le $\QQ$-espace vectoriel $\Gamma$ d'un ordre quelconque compatible avec sa loi de groupe, qui permet de considérer les termes de plus haut degré des expressions qu'on manipule. 

\trois{idmaxgrad} Un idéal homogène $ I$ d'un anneau gradué $ A$ est dit premier (resp. maximal) s'il est différent de $ A$ et si l'on a $$ab\in  I\Rightarrow a\in  I\;{\rm ou} \;b\in  I$$ pour tout couple $(a,b)$ d'éléments homogènes de $ A$ (resp. et s'il est maximal parmi les idéaux homogènes stricts de $ A$). L'idéal $ I$ est premier (resp. maximal) si et seulement si l'anneau gradué $ A/ I$ est intègre (resp. est un corps gradué). Tout anneau gradué non nul admet un idéal homogène maximal. 

\medskip
La remarque \ref{remcorpsgrad} ci-dessus se décline dans ce contexte : un idéal homogène de $ A$ est premier si et seulement si il est premier comme idéal ; il peut par contre être maximal en tant qu'idéal homogène, sans l'être en tant qu'idéal tout court. 

\trois{principagrad} Un élément homogène d'un anneau gradué intègre est dit {\em irréductible} s'il est non inversible et s'il ne peut s'écrire comme produit non trivial de deux éléments homogènes (on vérifie comme au \ref{remcorpsgrad} qu'il ne peut dans ce cas pas non plus s'écrire comme produit non trivial de deux éléments quelconques). 

\medskip
Un anneau gradué $ A$ est dit {\em principal} s'il est intègre et si tout idéal homogène $ I$ de $ A$ est de la forme $(a)$ pour un certain élément homogène $a$ de $ A$ (qui est alors uniquement déterminé à un inversible près). Dans un anneau gradué principal tout élément homogène non nul s'écrit de manière unique comme produit d'éléments homogènes irréductibles, les notions de PGCD et PPCM de deux éléments homogènes sont bien définies, les idéaux homogènes premiers non nuls sont maximaux, un élément homogène non nul est irréductible si et seulement si l'idéal homogène qu'il engendre est premier, etc.

\deux{anneaugrad} Soit $ A$ un anneau gradué et soit $(r_{1},\ldots,r_{n})\in \Gamma^{n}$. On note $$ A[r_{1}^{-1}T_{1},\ldots,r_{n}^{-1}T_{n}]$$ (en omettant éventuellement les $r_{i}$ égaux à $1$) l'anneau gradué dont l'anneau sous-jacent est $ A[T_{1},\ldots,T_{n}]$, et dont le groupe des éléments homogènes de degré $s$ est formé, pour $s$ fixé, des polynômes de la forme $\sum a_{I}{\bf T}^{I}$ avec $a_{I}$ homogène de degré $s{\bf r}^{-I}$ pour tout $I$. 

\medskip
Par abus, on notera $ A[[\tau_{1},\ldots,\tau_{r}]]$ l'anneau gradué $$\bigoplus_{s} A_{s}[[\tau_{1},\ldots,\tau_{r}]]$$ où $ A_{s}$ désigne pour tout $s$ l'ensemble des éléments de $ A$ de degré $s$. 

\deux{excorpsgrad} Soit $ k$ un corps gradué ({\em e.g.} un vrai corps muni de la graduation triviale), et soit $r\in \Gamma$ tel que $ k_\gamma=\{0\}$ pour tout $\gamma\in r^\ZZ$. L'anneau gradué $$ k[r^{-1}T, rT^{-1}]:= k[r^{-1}T, rS]/(TS-1)$$ est un corps gradué qui n'est jamais un corps (au sens usuel) : par exemple, l'élément $1+T$ n'y est pas inversible. 

\deux{polyprincgrad} Si $ k$ est un corps gradué et si $\gamma\in\Gamma$ l'anneau gradué $ k[T/\gamma]$ est principal. Si $P$ est un élément homogène de $ k[T/\gamma]$ son degré  désignera son degré relatif à la graduation de $ k[T/\gamma]$ ; lorsqu'on voudra évoquer son degré au sens usuel pour un polynôme, on parlera de degré {\em monomial}. 

\medskip
Si $P$ est un élément homogène et unitaire de $ k[T/\gamma]$ et si $n$ désigne le degré monomial de $P$ alors le degré de $P$ est égal à $\gamma^n$. 

\medskip
Si $P$ est un élément homogène de $ k[T/\gamma]$ une {\em racine} de $P$ dans $ k$ est un élément $x$ de $ k$ homogène de degré $\gamma$ tel que $P(x)=0$. Un élément $x$ de $ k$ homogène de degré $\gamma$ est racine de $P$ si et seulement si $T-x$ divise $P$. Si $P$ est non nul et si $n$ désigne son degré monomial, il a au plus $n$ racines dans $ k$ (comptées avec multiplicité) ; il en a exactement $n$ si et seulement si il s'écrit sous  la forme $a\prod (T-x_i)$ où $a$ et les $x_i$ sont des éléments homogènes de $ k$, les $x_i$ étant de surcroît tous de degré $\gamma$. On dit alors que $P$ est scindé dans $ k$. 

\deux{entiers} Soit $ A$ un anneau gradué. Une {\em $ A$-algèbre graduée} est un anneau gradué  $ B$ muni d'un morphisme $ A\to  B$ d'anneaux gradués. Si $ A$ est un corps gradué et si $ B$ est une $ A$-algèbre graduée non nulle alors $ A\to  B$ est injective ; si de plus $ B$ est un corps gradué, on dit que c'est une {\em extension graduée} de $ A$.

\deux{espvectgrad} {\em Modules et espaces vectoriels gradués.} Soit $ A$ un anneau gradué. Un {\em $ A$-module gradué} est un $ A$-module (classique) $ M$ muni d'une décomposition en somme directe $\bigoplus \limits_{\gamma \in \Gamma}  M_\gamma$ telle que $ A_\gamma. M_\delta\subset  M_{\gamma\delta}$ pour tout $(\gamma, \delta)$ ; on dispose d'une notion évidente d'application linéaire graduée de degré donné $\delta$ entre deux $ A$-modules gradués, de sous-module gradué, de somme directe de modules gradués, de produit tensoriel de modules gradués, etc. 

\medskip
Si $\gamma \in \Gamma$ on note $ A(\gamma)$ le $ A$-module gradué dont le $ A$-module sous-jacent est $ A$ et dont la graduation est telle que $ A(\delta)_\gamma= A_{\gamma\delta}$ pour tout $\delta$. Si $ M$ est un $ A$-module gradué et si $(m_i)$ est une famille d'éléments homogènes de $ M$ de degrés respectifs $\gamma_i$, il existe une unique application linéaire graduée et de degré $1$ de $\bigoplus  A(\gamma_i^{-1})$ dans $ M$ qui envoie $1_i$ sur $m_i$ pour tout $i$ ; on dit que la famille $(m_i)$ est libre (resp. génératrice, resp. une base) si cette application est injective (resp. surjective, resp. bijective). Un $ A$-module est dit de type fini s'il admet une famille génératrice finie. 

\medskip
Si $ k$ est un corps gradué, on parlera de $ k$-espace vectoriel gradué plutôt que de $ k$-module gradué. Si $ M$ est un $ k$-espace vectoriel gradué, il possède une base, et toutes ses bases ont même cardinal ; ce dernier est appelé la {\em dimension} de $ M$.

\deux{diffgrad} Si $ A$ est un anneau gradué, si $ B$ est une $ A$-algèbre graduée, et si $ M$ est un $ B$-module gradué une {\em $ A$-dérivation} de $ B$ dans $ M$ est une application $ A$-linéaire graduée ${\rm d} :  B\to  M$ de degré $1$ telle que ${\rm d}a=0$ pour tout élément homogène $a\in  A$ et $${\rm d}(bb')=b\;{\rm d}b'+b'{\rm d}b$$ pour tout couple $(b,b')$ d'éléments homogènes de $ B$. La catégorie des $ B$-modules munis d'une $ A$-dérivation de source $ B$ admet un objet initial que l'on note $\Omega_{ B/ A}$ ; c'est le {\em module gradué des différentielles de $ B$ sur $ A$}. 

\deux{nakagrad} Un anneau gradué $ A$ est dit {\em local} s'il possède un et un seul idéal homogène maximal. On dispose d'une variante graduée du lemme de Nakayama : si $ A$ est un anneau local gradué d'idéal maximal $\got m$ et si le $ A$-module gradué $\got m$ est de type fini, une famille finie $(m_i)$ d'éléments homogènes de $\got m$ engendre $\got m$ si et seulement si les classes des $m_i$ engendrent le $ k$-espace vectoriel gradué $\got m/\got m^2$, où $ k= A/\got m$. Le cardinal minimal d'une famille génératrice de $\got m$ est donc égal à la dimension de $\got m/\got m^2$. 

\deux{clotintegr} Soit $ A$ un anneau gradué et soit $ B$ une $ A$-algèbre graduée. Si $b$ est un élément de $ B$ homogène de degré $\gamma$, on dit que $b$ est {\em entier sur $ A$} s'il existe un entier $n$ et un élément $P\in  A[T/\gamma]_{\gamma^n}$ unitaire et de degré monomial égal à $n$ tel que $P(b)=0$ ; il revient au même ({\em via} par exemple un avatar gradué du théorème de Cayley-Hamilton) de demander que $b$ vive dans une sous-$ A$-algèbre graduée de $ B$ qui soit de type fini comme $ A$-module gradué. L'ensemble des éléments homogènes de $ B$ entiers sur $ A$ est l'ensemble des éléments homogènes d'une sous-$ A$-algèbre graduée de $ B$ que l'on appelle la {\em fermeture intégrale} de $ A$ dans $ B$. 

\medskip
Soit $ A$ un anneau gradué et soit $ B$ une $ A$-algèbre graduée que l'on suppose entière ({\em i.e.} constituée d'éléments entiers), soit $\got p$ un idéal premier homogène de $ A$ et soit $\got q$ un idéal premier homogène de $ B$ tel que $\got q\cap  A=\got p$. Sous ces hypothèses : 

\medskip
1) $\got p$ est maximal si et seulement si $\got q$ est maximal ; 

2) si $\got p'$ est un idéal premier de $ A$ contenant $\got p$, il existe un idéal premier $\got q'$ de $ B$ contenant $\got q$ tel que $\got q'\cap  B=\got p'$ (c'est le {\em going-up} gradué). 

\subsection*{La théorie de Galois graduée} 

\deux{extgrad} Soit $ k$ un corps gradué et soit $ L$ une extension de $ k$. Un élément homogène $x$ de $ L$ est dit {\em algébrique} sur $ k$ s'il est entier sur $k$ (au sens de \ref{clotintegr}). Si $x$ est de plus non nul et si $\gamma$ désigne son degré, l'unique générateur homogène et unitaire de l'annulateur de $x$ dans $ k[T/\gamma]$ est appelé le {\em polynôme minimal} de $x$ sur $ k$ ; il est irréductible ; son degré monomial est simplement appelé le {\em degré de $x$ sur $ k$.} 

\medskip
Une extension graduée $L$ de $k$ est dite {\em algébrique} si tous les éléments homogènes de $L$ sont algébriques sur $k$. Si $L$ est une extension graduée algébrique de $k$, le groupe des $k$-automorphismes de $L$ sera noté $\mathsf{Gal}(L/k)$.

\medskip
{\em Remarque.} Nous n'avons pas défini le polynôme minimal de $0$ pour la raison suivante : ce dernier devrait être égal à $T$, mais pourrait être vu comme appartenant à $ k[T/\gamma]$ pour n'importe quel $\gamma$ ; c'est cette ambiguïté qui nous a conduit à ne pas le considérer. 

\deux{ruptdecomp} Soit $ k$ un corps gradué et soit $\gamma\in \Gamma$. Si $P$ est un élément homogène et irréductible de $ k[T/\gamma]$ un {\em corps de rupture} de $P$ est une extension graduée $ L$ de $ k$ engendrée par une racine de $P$ ; si $P$ est un élément homogène non nul de $ k[T/\gamma]$, un {\em corps de décomposition} de $P$ est une extension graduée $ L$ de $ k$ dans laquelle $P$ est scindée, et qui est engendrée par les racines de $P$. 

\medskip
Si $P$ est un élément homogène et irréductible de $ k[T/\gamma]$, il admet un corps de rupture, à savoir $ k[T/\gamma]/P$ ; on en déduit que si $P$ est un élément homogène non nul de $ k[T/\gamma]$, il admet un corps de décomposition. Combiné avec un raisonnement fondé sur le lemme de Zorn, ceci entraîne l'existence d'une {\em clôture algébrique graduée} de $ k$, c'est-à-dire d'une extension algébrique graduée $ L$ de $ k$ telle que pour tout $\gamma\in \Gamma$, tout élément homogène et non nul  $P$ de $ L[T/\gamma]$ soit scindé dans $ L$. 

\deux{dedekgrad} Soit $ k$ un corps gradué, soit $ A$ une $ k$-algèbre graduée et soit $ L$ une extension graduée de $ k$. Si $\phi_1,\ldots, \phi_n$ sont des morphismes deux à deux distincts de $ A$ dans $ L$, ils forment une famille libre dans l'ensemble des applications $ k$-linéaires graduées (de degré quelconque) de $ A$ dans 
$ L$ (c'est la variante graduée du lemme d'indépendance des caractères de Dedekind). On en déduit, en appliquant ce lemme sur le corps gradué $ L$ à la famille des application $\psi_i :  L\otimes_{ k} A\to  L$ induites par les $\phi_i$, que si $ A$ est de dimension finie $d$ sur $ k$, il existe au plus $d$ morphismes de $ A$ dans $ L$.

\trois{borneplong} En particulier, si $ F$ est une extension graduée finie de $ k$ l'ensemble des $ k$-plongements de $ F$ dans $ L$ est fini, de cardinal majoré par $[ F: k]$. 

\trois{bornecorpsdcec} Supposons de plus que $ F$ soit un corps de décomposition d'un élément homogène $P$ non nul de $ k[T/\gamma]$ pour un certain $\gamma$. Dans ce cas, si $P$ est scindé (resp. scindé à racines simples) dans $ L$ il y a au moins un (resp. exactement $[ F: k]$) plongement(s) de $ F$ dans $ L$. Cela peut se démontrer par récurrence sur $[F:k]$, le point clef étant que si $Q$ est un élément homogène irréductible de $k[T/\gamma]$ alors $T\mapsto \phi(T)$ établit une bijection entre $\mathsf{Hom}_k(k[T/\gamma]/Q,L)$ et l'ensemble des racines de $Q$ dans $L$. 

\trois{automorphgrad} Il s'ensuit que deux corps de décomposition d'un même élément homogène $P$ de $ k[T/\gamma]$ sont isomorphes ; on en déduit, à l'aide du lemme de Zorn, que deux clôtures algébriques graduées de $ k$ sont isomorphes. 

\deux{algetale} {\em Algèbres étales graduées}. Soit $ k$ un corps gradué et soit $\overline { k}$ une clôture algébrique graduée de $ k$. 

\trois{polsep} Si $\gamma\in \Gamma$, un élément homogène non nul $P\in  k[T/\gamma]$ est dit {\em séparable} s'il est premier à son polynôme dérivé $P'$ (qui est homogène de degré $\delta/\gamma$ si $P$ est homogène de degré $\delta$) ; il revient au même de demander que les racines de $P$ dans $\overline { k}$ soient toutes simples, ou encore que $\overline { k}[T/\gamma]/P$ soit réduit. Si $P$ est irréductible, il est séparable si et seulement si $P'\neq 0$ ; en général, $P$ est séparable si et seulement si il n'a que des facteurs irréductibles simples, et si chacun d'eux est séparable. 

Soit $L$ une extension graduée algébrique de $k$. Soit  $x$ un élément homogène de $L$. On dit que $x$ est {\em séparable} sur $k$ s'il existe $\gamma\in \Gamma$ et un élément homogène et séparable $P$ de $k[T/\gamma]$ dont $x$ est une racine ; lorsque $x\neq 0$, cela revient à demander que son polynôme minimal soit séparable. 

\medskip
Si tous les éléments homogènes de $L$ sont séparables sur $k$, on dit que $L$ elle-même est séparable sur $k$. 

\trois{intralget} Si $ A$ est une $ k$-algèbre graduée finie ({\em i.e.} de dimension finie sur $ k$), elle est isomorphe à un produit fini de $ k$-algèbres graduées locales artiniennes ({\em i.e.} dont l'idéal homogène maximal est constitué d'éléments nilpotents). Si $ A$ est locale, si $\got m$ désigne son idéal homogène maximal et si $ A/\got m= k$ alors $\Omega_{ A/ k}$ s'identifie au dual gradué de $\got m/\got m^2$. 

\trois{algetkbarre} On en déduit à l'aide du lemme de Nakayama gradué que pour une $\overline{  k}$-algèbre graduée finie locale (resp. graduée finie quelconque) $A$ les assertions suivantes sont équivalentes : 

1) $ A$ est réduite ; 

2) $ A$ est égale à $\overline{  k}$ (resp. est le produit d'un nombre fini de copies de $\overline k$) ; 

3) $\Omega_{ A/\overline { k}}=0$.

\trois{algetgen} Si $ A$ est une $ k$-algèbre graduée finie, on dit que $ A$ est {\em étale} si $\Omega_{ A/ k}=0$. Par ce qui précède, il revient au même de demander que $\overline{  k}\otimes_{ k}  A$ soit réduite, ou encore que  $\overline { k}\otimes_{ k}  A$ soit un produit fini de copies de $\overline { k}$. 

\medskip
Indiquons quelques propriétés qui résultent immédiatement de la définition. 

Si $ A$ est une $ k$-algèbre graduée étale alors elle est réduite (puisqu'elle se plonge dans $\overline{ k}\otimes_{ k} A$) et est donc un produit de corps. 

Si $ A$ et $ B$ sont deux $ k$-algèbres graduées finies alors $ A\times  B$ est étale si et seulement si $ A$ et $ B$ sont étales.

Si $ A$ est une $ k$-algèbre graduée étale et si $ B$ est une sous-algèbre graduée de $ A$ alors $ B$ est étale : cela résulte du fait que $\overline{ k}\otimes_{ k} B$ est alors une sous-algèbre graduée de $\overline{ k}\otimes_{ k} B$, et est en particulier réduite.

\trois{algetsep} Il résulte de ce qui précède qu'une $ k$-algèbre graduée finie est étale si et seulement si elle s'écrit comme un produit d'extensions graduées finies de $ k$ qui sont étales. Pour comprendre ce que sont les $ k$-algèbres graduées étales, il suffit donc de comprendre à quelle condition une extension graduée finie $ L$ de $ k$ est étale. Soit donc $ L$ une extension graduée finie de $ k$. 

\medskip
Supposons que $ L$ soit étale, soit $x$ un élément homogène non nul de $ L$ dont on note $\gamma$ le degré et soit $P$ son polynôme minimal. En tant que sous-algèbre graduée de $L$, la $ k$-algèbre graduée $ k[x]\simeq  k[T/\gamma]/P$ est étale, c'est-à-dire qu'elle reste réduite après extension des scalaires à $\overline{ k}$  ; par conséquent, $P$ est séparable ; on en conclut que $ L$ est séparable. 

\medskip
Réciproquement, soit $ L$ une extension graduée finie de $ k$ engendrée comme $ k$-algèbre graduée par une famille finie $(x_1,\ldots, x_n)$ d'éléments homogènes que l'on suppose séparables sur $k$ ; quitte à éliminer les $x_i$ égaux à $0$, on peut faire l'hypothèse qu'ils sont tous non nuls. Pour tout $i$, notons $P_i$ le polynôme minimal de $x_i$. Comme $x_i$ est séparable, $P'_i$ est premier à $P_i$ et $P'_i(x_i)$ est donc non nul ; l'égalité $P_i(x_i)=0$ entraîne que $P'_i(x_i){\rm d}x_i=0$ et donc que ${\rm d}x_i=0$ ; ceci valant pour tout $i$, on a $\Omega_{ L/ k}=0$, et $ L$ est étale. 

\medskip
Ainsi, une $ k$-algèbre graduée finie est étale si et seulement si elle s'écrit comme un produit fini d'extensions graduées séparables de $ k$. Et pour qu'une extension graduée finie $ L$ de $ k$ soit séparable, il suffit qu'elle soit engendrée par une famille finie d'éléments séparables.

\trois{fermesp} Énonçons maintenant quelques conséquences de ce qui précède. Si $ L$ est une extension graduée finie séparable de $ k$ et si $ F$ est une extension graduée finie séparable de $ L$ alors $ F$ est une extension graduée finie séparable de $ k$. 

Si $ L$ est une extension graduée quelconque de $ k$, l'ensemble des éléments homogènes de $ L$ qui sont séparables sur $ k$ est l'ensemble des éléments homogènes d'une extension graduée algébrique de $ k$ appelée {\em fermeture séparable} de $ k$ dans $ L$ ; celle-ci ne possède aucune extension séparable stricte dans $ L$.

Si $ L$ est une extension algébrique de $ k$, la fermeture séparable de $ k$ dans $ L$ est réduite à $ k$ si et seulement si $ L$ est une extension {\em radicielle} de $ k$, c'est-à-dire si et seulement si le polynôme minimal de tout élément homogène non nul de $ L$ est de la forme $T^{p^n}-a$ où $a$ est un élément homogène non nul de $ k$ et où $p$ est l'exposant caractéristique de $ k$.

\medskip
La fermeture séparable $ k^s$ de $ k$ dans $\overline { k}$ est une extension graduée séparable de $ k$ n'admettant aucune extension graduée finie séparable stricte. Une telle extension est appelée une {\em clôture séparable de $ k$} ; deux clôtures séparables de $ k$ sont isomorphes. 

\deux{autfixe} Soit $ L$ une extension graduée de $ k$ et soit $\mathsf G$ un groupe fini d'automorphismes de $ L$ ; supposons que $L^{\mathsf G}=k$. Soit $x$ un élément homogène non nul de $L$ et soit $\gamma$ son degré. Soit $P$ l'élément $\prod\limits_{y\in \mathsf G.x} (T-y)$ de $L[T/\gamma]$, qui est homogène de degré $\gamma^{\sharp\mathsf G.x}$. Le polynôme $P$ est invariant sous $\mathsf G$ et appartient donc à $k[T/\gamma]$ ; il est scindé à racines simples dans $L$ et annule $x$ ; par conséquent, $x$ est séparable sur $k$, et $\mathsf G.x$ engendre une sous-extension graduée de $L$ finie et séparable sur $k$ qui est stable sous $\mathsf G$. 

Ainsi, $L$ est réunion de ses sous-extensions graduées finies et séparables sur $k$ stables sous $\mathsf G$. Soit $F$ l'une d'elle et soit $n$ sa dimension sur $k$. Le groupe $\mathsf G$ agit sur $\overline k \otimes_k F\simeq \overline k^n$ ; il découle de l'égalité $L^{\mathsf G}=k$ que l'ensemble des éléments de $\overline k \otimes_k F\simeq \overline k^n$ invariants sous $\mathsf G$ est la diagonale $\overline k.(1,\ldots, 1)$. 

Par ailleurs, le groupe des automorphismes de la $\overline k$-algèbre graduée $\overline k^n$ s'identifie à $\got S_n$, agissant par permutation des facteurs : cela résulte (par exemple) du fait que les éléments de la forme $(0,\ldots,0,1,0,\ldots, 0)$ de $\overline k^n$ sont ses idempotents homogènes non nuls minimaux, et sont donc permutés par tout automorphisme. Comme la sous-algèbre de $\overline k^n$ formée des éléments invariants sous $\mathsf G$ est égale à la diagonale, $\mathsf G$ agit transitivement sur $\{1,\ldots, n\}$, ce qui oblige $n$ à être inférieur ou égal au cardinal de $\mathsf G$ ; ainsi, $[F:k]\leq \sharp \mathsf G$. 

\medskip
Par conséquent, $L$ est une extension graduée finie séparable de $k$, et l'on a la majoration $[L:k]\leq \sharp\mathsf G$. Comme par ailleurs $\mathsf G$ se  plonge dans le groupe des $k$-automorphismes de $L$, lequel a un cardinal borné par $[L:k]$, on en déduit que $[L:k]=\sharp\mathsf G$ et que $\mathsf G$ est égal au groupe des $k$-automorphismes de $L$. 

\deux{thogalgrad} Les ingrédients sont désormais réunis pour disposer d'une théorie de Galois graduée, qui se déduit formellement de tout ce qui précède. Nous allons en donner une formulation classique, et une plus catégorique, fidèle au point de vue de Grothendieck. 

\trois{theogalgradclass} {\em La formulation classique.} Soit $k$ un corps gradué et soit $L$ une extension graduée de $k$. Les assertions suivantes sont équivalentes : 

\medskip
i) il existe un groupe fini $\mathsf G$ d'automorphismes de $L$ tel que $k=L^{\mathsf G}$ ; 

ii) l'extension graduée $L/k$ est finie, séparable, et normale (ce qui signifie que si $x$ est un élément homogène non nul de $L$, son polynôme minimal $P$ est scindé dans $L$) ; 

iii) il existe $\gamma\in \Gamma$ et un élément homogène, non nul et séparable $P\in k[T/\gamma]$ tel que $L$ soit un corps de décomposition de $P$ sur $k$.

\medskip
De plus si i) est vraie alors $\mathsf G=\mathsf {Gal}(L/k)$ et $[L:k]=\sharp \mathsf G$, et si ii) ou iii) est vraie alors $k=L^{\mathsf{Gal}(L/k)}$. 

\medskip
Lorsqu'une extension graduée $L$ de $k$ satisfait ces conditions équivalentes, on dit qu'elle est {\em finie galoisienne}. Si $L$ est une extension graduée finie galoisienne de $k$, on dispose sur $L$ de la {\em correspondance de Galois} : les applications $\mathsf H\mapsto L^{\mathsf H}$ et $F\mapsto \mathsf{Gal}(L/F)$ établissent une bijection décroissante entre l'ensemble des sous-groupes de $\mathsf{Gal}(L/k)$ et celui des sous-extensions graduées de $L/k$. Si $\mathsf H$ est un  sous-groupe de $\mathsf{Gal}(L/k)$ et si $F$ est la sous-extension qui lui correspond alors $F$ est une extension galoisienne de $k$ si et seulement si $\mathsf H$ est un sous-groupe distingué de $\mathsf{Gal}(L/k)$, ce qui revient à demander que $F$ soit stable sous $\mathsf{Gal}(L/k)$ ; si c'est le cas, la restriction des automorphismes induit un isomorphisme $\mathsf{Gal}(L/k)/\mathsf H\simeq \mathsf{Gal}(F/k)$. 

\trois{theogalcat} {\em La formulation catégorique.} Soit $\mathsf G$ le groupe $\mathsf {Gal}(k^s/k)$ ; c'est un groupe profini qui coïncide avec $\mathsf{Gal}\; (\overline k/k)$ (cela provient du fait qu'une extension graduée radicielle n'a pas d'automorphismes non triviaux). 

\medskip
Le foncteur $A\mapsto \mathsf{Hom}_k(A,k^s)$ établit une équivalence entre la catégorie des $k$-algèbres graduées étales et celle des $\mathsf G$-ensembles discrets finis (un $\mathsf G$-ensemble discret est un ensemble muni d'une action continue de $\mathsf G$ pour laquelle les stabilisateurs sont ouverts). 

\deux{remelprem}  {\em Remarque.} Ainsi, la théorie de Galois classique survit à peu près mot pour mot dans le contexte gradué. Mentionnons toutefois un théorème qui n'est plus vrai : celui de l'élément primitif  ; nous allons en effet construire un contre-exemple à celui-ci. 

\medskip
Soit $k_0$ un vrai corps, que l'on gradue trivialement, et supposons donnée une famille libre $(\gamma,\delta)$ du $\QQ$-espace vectoriel $\Gamma$. Désignons par $k$ le corps gradué $k_0[T/\gamma, \gamma T^{-1}, S/\delta, \delta S^{-1}]$. Soit $n$ un entier strictement supérieur à $1$ et inversible dans $k_0$, et soit $L$ l'algèbre graduée $k[U/\gamma^{1/n}, V/\delta^{1/n}]/(U^n-T, V^n-S)$.

On vérifie aisément que $L$ est une extension graduée finie et séparable de degré $n^2$ de $k$. Soit $x$ un élément homogène de $L$ ; il est de la forme $aU^mV^{m'}$ où $a$ est un élément homogène de $k$ ; par conséquent, $x^n\in k$ et $x$ ne peut engendrer $L$ qui est de degré $n^2$ sur $k$.

\section{Réduction graduée et théorie de la ramification modérée} 

On désigne maintenant par $\Gamma$ un groupe abélien divisible, sans torsion et {\em ordonné}, toujours noté multiplicativement. On lui adjoint un élément $0$ absorbant pour la multiplication. En pratique, $\Gamma$ sera égal à $\RR^*_+$. 

\medskip
On fixe un corps $k$ muni d'une valeur absolue ultramétrique $|.|$ prenant ses valeurs dans $\Gamma\cup\{0\}$, c'est-à-dire une application $|.|: k\to \Gamma\cup\{0\}$ possédant les propriétés suivantes : 

\medskip
$\bullet$ $|x|=0\iff x=0$ ; 

$\bullet$ $|1|=1$ ; 

$\bullet$ $|ab|=|a|.|b|$ pour tout $(a,b)$ ; 

$\bullet$ $|a+b|\leq \max (|a|,|b|)$ pour tout $(a,b)$. 

\medskip
On vérifie aussitôt que $|a+b|=\max(|a|, |b|)$ dès que $|a|$ et $|b|$ diffèrent. 

\medskip
\deux{semnorm} Si  $A$ est une $k$-algèbre, une {\em semi-norme} sur $A$ est une application $||.||$ de $A$ dans $\Gamma\cup\{0\}$ telle que $||a+b||\leq \max(||a||, ||b||)$ pour tout $(a,b)$ et telle que $||\lambda a||=|\lambda|.||a||$ pour tout $\lambda\in k$ et tout $a\in A$ ; là encore, on a $||a+b||=\max(||a||, ||b||)$ dès que $||a||$ et $||b||$ diffèrent. Si de plus $||a||=0\iff a=0$, on dira que $||.||$ est une {\em norme}. 

\deux{ineg} Soit $A$ une $k$-algèbre munie d'une semi-norme $||.||$. Pour tout $r\in \Gamma$, on note $A_{\leq r}$ (resp. $A_{<r}$) le sous-groupe de $A$ formé des éléments $a$ tels que $||a||\leq r$ (resp. $||a||<r$). 

On désigne par $\red A$ la {\em réduction graduée} de $A$ au sens de Temkin, c'est-à-dire l'algèbre $\Gamma$-graduée $\bigoplus\limits_{r\in \Gamma}A_{\leq r}A_{<r}$. Si $a\in A$ et si $r$ est un élément de $\Gamma$ tel que $||a||\leq r$, on notera $\red{a}_{r}$ l'image de $a$ dans $\red{A}_{r}$ ; si $r=||a||$, on écrira simplement $\red{a}$ ; si $||a||=0$, on pose $\red{a}=0$. 

\medskip
Remarquons que $\red A_1$ n'est autre que la réduction classique de $A$. 

\medskip
Lorsqu'on applique cette construction à $A=k$, l'anneau gradué $\red k$ obtenu est un corps gradué, appelé {\em corps gradué résiduel} de $k$ ; le corps résiduel de $k$ au sens classique est égal à $\red k_1$. On désigne par $p$ l'exposant caractéristique de $\red k$ (qui est égal à celui de $\red k_1$). 

\medskip
\deux{libregredk} {\bf Lemme.} {\em Soit $(a_{i})$ une famille d'éléments de $A$. Les propositions suivantes sont équivalentes :

\medskip
$i)$ $(\red{a_{i}})$ est une famille libre du $\red k$-espace vectoriel gradué $\red A$ ;

\medskip
$ii)$ Les normes des $a_{i}$ sont toutes non nulles, et $||\sum \lambda_{i}a_{i}||=\max |\lambda_{i}|.||a_{i}||$ pour toute famille $(\lambda_{i})$ d'éléments de $k$.~$\Box$ }

\medskip
\deux{dimgreda} {\bf Corollaire.} {\em On a l'inégalité $\dim{\red{k}}\red A\leq \dim{k}A$.} 

\medskip
{\em Démonstration.} Cela résulte du fait que toute famille d'éléments de $A$ qui vérifie $ii)$ est libre.~$\Box$

\subsection*{Corps résiduels gradués et théorie de Galois}

\deux{hypkhens} On suppose maintenant que le corps valué $k$ est hensélien, dans le sens suivant : {\em pour toute extension finie $L$ de $k$, il existe une unique valeur absolue sur  $L$ à valeurs dans $\Gamma\cup\{0\}$ et qui prolonge celle de $k$.} Si $L$ est une extension finie de $k$, on la considèrera implicitement comme munie de l'unique prolongement de $|.|$, que l'on notera encore $|.|$. Si $E$ est une $k$-algèbre étale, elle est isomorphe à un produit fini $\prod L_i$, où $L_i$ est pour tout $i$ une extension séparable de $k$, et nous la considérerons comme une algèbre normée {\em via} l'application $$(x_1,\ldots,x_n)\mapsto \max |x_i|.$$

\medskip
Commençons par établir différentes variantes graduées du lemme de Hensel.

\medskip
\deux{rhop} Soit $P=X^{n}+\sum_{i\leq n-1} a_{i}X^{i}$ un polynôme unitaire à coefficients dans $k$ ; soit $L$ un corps de décomposition de $P$. Posons $\rho(P)=\max |a_{i}|^{1/n-i}$. Un calcul immédiat montre que si $x\in L$ est tel que $|x|>\rho(P)$, alors $|x|^{n}
>|a_{i}|.|x|^{i}$ pour tout $i\leq n-1$ ; en conséquence, $x$ ne peut être racine de $P$. Par ailleurs, si $\rho(P)$ majorait strictement en valeur absolue toutes les racines de $P$ dans $L$, l'expression des $a_{i}$ en fonction de ces dernières fournirait l'inégalité absurde $\max |a_{i}|^{1/n-i}<\rho(P)$ ; il en découle que $\rho(P)$ coïncide avec le maximum des valeurs absolues des racines de $P$ dans $L$. Si $Q$ divise $P$, alors $\rho(Q)\leq \rho(P)$.

\medskip
Pour tout $r\in \Gamma$ majorant $\rho(P)$, on notera $\red{P}_{r}$ le polynôme $X^{n}+\sum \widetilde{a_{i,r^{n-i}}}X^{i}$ à coefficients dans $\red{k}$ ; c'est un élément de $\red{k}[X/r]$ homogène de degré $r^{n}$. Si $P$ s'écrit $RS$ avec $R$ et $S$ unitaires, on a a alors $\red{P}_{r}=\red{R}_{r}\red{S}_{r}$ ; en particulier, si l'on écrit $P=\prod (X-x_{i})$ dans $L$, alors $\red{P}_{r}=\prod (X-\widetilde{x_{i,r}})$. On en déduit que $r>\rho(P)$ si et seulement si $\red{P}_{r}=X^{n}$.

\medskip
Si $\rho(P)>0$ (autrement dit si $P$ n'est pas une puissance de $X$), on écrira simplement $\red{P}$ au lieu de $\red{P}_{\rho(P)}$.  

\medskip
Si $r\in \Gamma$ et si $R$ est un élément unitaire et homogène de $\red{k}[X/r]$, on appellera {\em relèvement admissible} de $R$ tout polynôme unitaire $\cal R$ appartenant à $k[X]$, tel que $\rho({\cal R})\leq r$ et tel que $\red{\cal R}_{r}=R$ ; l'existence d'un tel relèvement est immédiate. 

\medskip
\deux{prehensel} {\bf Lemme.} {\em Soit $P$ un polynôme unitaire à coefficients dans $k$ et soit $r$ un élément de $\Gamma$ supérieur ou égal à $\rho(P)$. Supposons que $\red{P}_{r}$ admet une factorisation $\red{P}=RS$, où $R$ et $S$ sont des éléments homogènes, unitaires, non constants et {\em premiers entre eux} de $\red{k}[X/r]$ ; sous cette hypothèse $\rho(P)=r$ et $P$ n'est pas irréductible sur $k$.} 

\medskip
{\em Démonstration.} Comme $R$ et $S$ sont premiers entre eux, $\red{P}_{r}$ est différent de $X^{n}$ et l'on a donc $\rho(P)=r$. Il reste à montrer que $P$ n'est pas irréductible ; on raisonne par l'absurde en supposant qu'il l'est. 

\medskip
Soit $L$ un corps de décomposition de $P$ ; écrivons $P=\prod (X-x_{i})$ dans $L$. Puisque $\prod(X-\widetilde{x_{i,r}})=RS$, il existe deux indices $i$ et $j$ distincts tels que $\widetilde{x_{i,r}}$ (resp. $\widetilde{x_{j,r}}$) soit une racine de $R$ (resp. de $S$) ; comme $R$ et $S$ sont premiers entre eux, $\widetilde{x_{j,r}}$ n'annule pas $R$. 

\medskip
Soit $\cal R$ un relevé admissible de $R$. Si $m$ désigne de degré de $R$ en $X$, alors $|{\cal R}(x_{i})|<r^{m}$ puisque $R(\widetilde{x_{i,r}})=0$. 

\medskip
Soit $g$ un $k$-automorphisme de $L$ envoyant $x_{i}$ sur $x_{j}$ ; on a nécessairement $|{\cal R}(x_{j})|=|g({\cal R}(x_{i}))|=|{\cal R}(x_{i})|<r^{m}$, ce qui est contradictoire avec le fait que $R(\widetilde{x_{j,r}})\neq 0$.~$\Box$ 

\medskip
\deux{hensel} {\bf Lemme.} {\em Soit $P$ un polynôme unitaire à coefficients dans $k$ et soit $r$ un élément de $\Gamma$ majorant $\rho(P)$. Supposons que $\red{P}_{r} $ s'écrit $\prod P_{i}$, où les $P_{i}$ sont des éléments homogènes, unitaires, non constants et {\em deux à deux premiers entre eux} de $\red{k}[X/r]$. Le polynôme $P$ admet alors une factorisation $P=\prod {\cal P}_{i}$ où chaque ${\cal P}_{i}$ est unitaire et tel que $\widetilde{{\cal P}_{i,r}}=P_{i}$.}

\medskip
{\em Démonstration.} On raisonne par récurrence sur le degré de $P$. S'il est nul le résultat est trivial ; on suppose donc que $\mbox{deg}\;P>0$ et que le résultat a été prouvé en degrés strictement inférieurs à $\mbox{deg}\;P$. Si le nombre de facteurs $P_i$ est égal à $1$, il n'y a rien à démontrer ; sinon, il résulte du lemme \ref{prehensel} que $P$ n'est pas irréductible ; écrivons alors $P=QR$, où $Q$ et $R$ sont unitaires et non constants ; notons que $\red{P}_{r}=\red Q_r\red R_r$.

\medskip
Pour tout $i$ notons $Q_i$ (resp. $R_i$) le PGCD unitaire de $\red Q$ (resp. $\red R$) et $P_i$. Les $Q_i$ (resp. $R_i$) sont homogènes, unitaires et deux à deux premiers entre eux ; on a $\red Q_r=\prod Q_i$, $\red R_r=\prod R_i$ et $Q_iR_i=P_i$ pour tout $i$. Comme $R$ et $Q$ sont tous deux de degré strictement inférieur à $\mbox{deg}\;P$, l'on peut leur appliquer l'hypothèse de récurrence ; on peut donc écrire $Q=\prod{\cal Q}_i$ et $R=\prod{\cal R}_i$ où ${\cal Q}_i$ (resp. ${\cal R}_i$) est pour tout $i$ un polynôme unitaire tel que $\red{{\cal Q}}_i=Q_i$ (resp.  $\red{{\cal R}}_i=R_i$). La famille des polynômes ${\cal P}_i:={\cal Q}_i{\cal R}_i$ convient alors.~$\Box$ 

\medskip
\deux{corohensel} {\bf Corollaire (lemme de Hensel gradué).} {\em  Soit $P$ un polynôme unitaire à coefficients dans $k$ et soit $r$ un élément de $\Gamma$ majorant $\rho(P)$. Si $\red{P}_{r}$ admet une racine simple $\lambda$ dans $\red{k}_{r}$, il existe alors une et une seule racine $l$ de $P$ dans $k$ telle que $\red{l}_{r}=\lambda.$~$\Box$ }

\medskip
\deux{ef} {\bf Proposition.} {\em Soit $L$ une extension finie de $k$ ; posons $e=[|L^{*}|:|k^{*}|]$ et $f=[\red{L}_1:\red{k}_1]$. On a alors $[L:k]\geq [\red{L}:\red{k}]=ef$. De plus, les conditions suivantes sont équivalentes : 

\medskip
$i)$ $\red{L}$ est une extension séparable de $\red{k}$ ;

\medskip
$ii)$ $\red{L}_1$ est une extension séparable de $\red{k}_1$ et $[|L^{*}|:|k^{*}|]$ est premier à $p$.}

\medskip
{\em Démonstration.} Pour tout élément $\gamma$ de $|k^{*}|$, choisissons $\lambda_{\gamma}$ dans $k^{*}$ tel que $|\lambda_{\gamma}|=\gamma$. Pour tout élément $\delta$ de $|L^{*}|/|k^{*}|$ donnons-nous $\mu_{\delta}\in L^{*}$ telle que l'image de $|\mu_{\delta}|$ dans $|L^{*}|/|k^{*}|$ soit égale à $\delta$. Fixons enfin une base $(\alpha_{1},\ldots,\alpha_{f})$ de $\red{L}_1$ sur $\red{k}_1$. Des deux décompositions $$\red{k}=\bigoplus_{\gamma}\red{k}_1.\red{\lambda_{\gamma}} \;\mbox{et}\;\red{L}=\bigoplus_{\gamma,\delta} \red{L}_1.\red{\mu_{\delta}}.\red{\lambda_{\gamma}},$$ on déduit que $(\alpha_{i}\red{\mu_{\delta}})$ est une base de $\red{L}$ sur $\red{k}$, d'où l'égalité $[\red{L}:\red{k}]=ef$. La majoration $[\red{L}:\red{k}]\leq [L:k]$ est fournie par le corollaire~\ref{dimgreda}.

\medskip
Supposons que $i)$ soit vraie ; le polynôme minimal sur $\red{k}$ de tout élément homogène non nul de $\red{L}$ est alors séparable ; c'est en particulier vrai pour les éléments de de degré $1$, c'est-à-dire ceux qui appartiennent à $\red{L}_1^{*}$ ; ainsi, $\red{L}_1$ est séparable sur $\red{k}_1$. Par ailleurs, soit $x$ un élément de $L^{*}$ tel que $|x|^{p}$ appartienne à $|k^{*}|$. Il existe $l\in\red{L}_1$ et $\omega\in \red{k}$ tel que $\red{x}^{p}=l\omega$ ; on en déduit que $\red{x}$ est purement inséparable sur le sous-corps gradué $\red{L}_1.\red{k}$ de $\red{L}$. Le corps $\red{L}$ est par hypothèse séparable sur $\red{k}$, il l'est {\em a fortiori} sur $\red{L}_1.\red{k}$. En conséquence, $\red{x}\in \red{L}_1.\red{k}$ ; ceci implique que le degré de $\red{x}$, autrement dit la valeur absolue de $x$, appartient à $|k^{*}|$ ; on a donc établi $ii)$. 

\medskip
Supposons que $ii)$ soit vraie. Le corps sous-corps gradué $\red{L}_1.\red{k}$ de $\red{L}$ est alors séparable sur $\red{k}$. Soit $\xi$ un élément homogène non nul de $\red{L}$ et soit $x$ un élément de $L^{*}$ tel que $\red{x}=\xi$. Soit $m$ un entier premier à $p$ tel que $|x|^{m}\in |k^{*}|$. Il existe $l\in \red{L}_1$ et $\omega\in \red{k}$ tel que $\xi^{m}=l\omega$ ; en conséquence, $\xi$ est séparable sur $\red{L}_1.\red{k}$. Le corps $\red{L}$ est donc séparable sur $\red{L}_1.\red{k}$, lequel est lui-même séparable sur $\red{k}$ ; dès lors, $\red{L}$ est séparable sur $\red{k}$.~$\Box$

\subsection*{La ramification modérée}

\medskip
On se propose maintenant de {\em reformuler} la théorie de la ramification modérée dans le langage des réductions graduées ; comme nous allons le voir, les extensions modérément ramifiés se comportent, dans ce cadre, comme les extensions non ramifiées dans le contexte classique -- on retrouve ainsi un phénomène bien connu en géométrie logarithmique, à laquelle le point de vue gradué est étroitement apparenté. 

\medskip
Ce qui suit est donc simplement une présentation nouvelle de résultats très classiques ; mais nous avons choisi, pour la commodité du lecteur, d'en redonner les démonstrations, que nous avons traduites en langage gradué pour qu'elles fournissent directement les énoncés requis.

\medskip
\deux{redquasigal} {\bf Proposition.} {\em Soit $L$ une extension finie galoisienne de $k$. 

\medskip
$a)$ Les extensions $\red{L}/\red{k}$ et $\red{L}_1/\red{k}_1$ sont normales, et les deux flèches naturelles $\mathsf {Gal}(L/k)\to \mathsf{Gal}(\red{L}/\red{k})$ et $\mathsf{Gal}(\red{L}/\red{k})\to\mathsf{Gal}(\red{L}_1/\red{k}_1)$ sont surjectives ;

\medskip
$b)$ La seconde flèche mentionnée ci-dessus s'insère dans une suite exacte $$1\to\mathsf{Hom}(|L^{*}|/|k^{*}|,\red{L}_1^{*})\to\mathsf {Gal}(\red{L}/\red{k})\to\mathsf {Gal}(\red{L}_1/\red{k}_1)\to 1. $$} 

\medskip
{\em Démonstration.} Prouvons $a)$. Soit $\xi$ un élément homogène non nul de $\red{L}$ ; relevons-le en un élément $x$ de $L$, de valeur absolue égale au degré $r$ de $\xi$. Soit $P$ le polynôme minimal de $x$ sur $L$. Écrivons $P=\prod (X-x_{i})$, où $x_{1},\ldots,x_{d}$ sont les éléments de l'orbite de $x$ sous $\mathsf {Gal}(L/k)$. Remarquons que $\rho(P)=r>0$. Le polynôme $\red{P}$ est donc bien défini, il est égal à $\prod (X-\red{x_{i}})$ ; c'est un élément unitaire homogène de $\red{k}[X/r]$ qui annule $\xi$ et est scindé dans $\red{L}$ ; le polynôme minimal $R$ de $\xi$ est donc scindé dans $\red{L}$, ce qui montre que $\red{L}$ est une extension normale de $\red k$. En se limitant au cas où $\xi$ est de degré $1$, on établit en particulier le caractère normal de $\red{L}_1/\red{k}_1$. 

\medskip
Soit $\eta$ une racine de $R$ distincte de $\xi$ ; elle coïncide avec $\red{x_{i}}$ pour un certain $i$. Il existe $g\in \mathsf {Gal}(L/k)$ tel que $g(x)=x_{i}$, et l'image de $g$ dans $\mathsf{Gal} (\red{L}/\red{k})$ envoie $\xi$ sur $\eta$. On en déduit qu'un élément homogène non nul de $\red{L}$ est invariant sous l'image de $\mathsf {Gal}(L/k)$ si et seulement si il coïncide avec tous ses conjugués, donc si et seulement si il est invariant sous l'image de $\mathsf{Gal} (\red{L}/\red{k})$ ; par la correspondance de Galois, l'image de $\mathsf {Gal}(L/k)$ est égale à $\mathsf{Gal} (\red{L}/\red{k})$ tout entier. Là encore, il suffit de se restreindre au cas où $\xi$ est de degré $1$ pour obtenir la surjectivité de $\mathsf {Gal}(L/k)\to \mathsf{Gal}(\red{L}/_1\red{k}_1)$.

\medskip
Établissons maintenant $b)$. La $|L^*|$-graduation de $\red L$ induit sur ce dernier une $|L^*|/|k^*|$-graduation $\red L=\bigoplus\limits_{\delta\in |L^*|/|k^*|} E_\delta$ ; pour tout $\delta$ l'ensemble $E_\delta-\{0\}$ est un $(\red L_1.\red k)^*$-torseur, le terme de torseur étant ici à prendre dans une acception homogène que nous laissons au lecteur le soin de préciser ; notons que $E_1=\red L_1.\red k$.

\medskip
Soit $\psi$ un morphisme de groupes de $|L^*|/|k^*|$ dans $\red L_1^*$. L'application de $\red L$ dans lui-même qui envoie $\sum x_\delta$ sur $\sum \psi(\delta)x_\delta$ est un $\red L_1.\red k$-automorphisme de $\red L$ ; la flèche $\mathsf{Hom}(|L^{*}|/|k^{*}|,\red{L}_1^{*})\to\mathsf {Gal}(\red{L}/\red L_1.\red{k})$ ainsi définie est injective par construction. 

\medskip
Montrons qu'elle est surjective. Soit $\phi\in \mathsf {Gal}(\red{L}/\red L_1.\red{k})$ et soit $\delta \in|L^*|/|k^*|$. Si $x$ et $y$ sont deux éléments non nuls et homogènes de $E_\delta$, leur quotient appartient à $\red L_1.\red k$ et l'on a donc $\phi(x)/\phi(y)=x/y$ ; ceci implique l'existence d'un élément homogène non nul $\psi(\delta)$ de $\red L_1.\red k$ tel que $\phi(x)=\psi(\delta)x$ pour tout $x\in E_\delta$ ; cette égalité appliquée à n'importe quel élément $x$ homogène et non nul de $E_\delta$ force le degré de $\psi(\delta)$ à être égal à $1$, et l'on a donc $\psi(\delta)\in \red L_1^*$. Du fait que $\phi$ est un morphisme on déduit immédiatement que $\psi\in\mathsf{Hom}(|L^{*}|/|k^{*}|,\red{L}_1^{*})$, d'où la surjectivité souhaitée.  

\medskip
Le groupe $\mathsf{Hom}(|L^{*}|/|k^{*}|,\red{L}_1^{*})$ s'identifie donc à $\mathsf {Gal}(\red{L}/\red L_1.\red{k})$, ce qu'il fallait démontrer.~$\Box$ 

\medskip
\deux{grouperam} Si $L$ est un extension finie galoisienne de $k$, on note $\mathsf{I}_{\mathsf {grad}}(L/k)$ le {\em groupe d'inertie gradué} de $L/k$, c'est-à-dire le noyau de $\mathsf {Gal}(L/k)\to \mathsf{Gal}(\red{L}/\red{k})$. Remarquons que par définition, $\mathsf{I}_{\mathsf {grad}}(L/k)$ est l'ensemble des $k$-automorphismes $g$ de $L$ tels que pour tout $x$ non nul dans $L$, l'on ait l'inégalité $|g(x)-x|<|x|$. 

\medskip
\deux{wcompatible} {\bf Lemme.} {\em Soit $L$ une extension finie galoisienne de $k$ et soit $F$ une sous-extension galoisienne de $L$. La suite $$1\to \mathsf{I}_{\mathsf {grad}}(L/F)\to\mathsf{I}_{\mathsf {grad}}(L/k)\to\mathsf{I}_{\mathsf {grad}}(F/k)\to 1$$ est exacte.} 

\medskip
{\em Démonstration.} Le seul fait non immédiat est la surjectivité de la flèche $\mathsf{I}_{\mathsf {grad}}(L/k)\to\mathsf{I}_{\mathsf {grad}}(F/k)$. Soit $g\in \mathsf{I}_{\mathsf {grad}}(F/k)$ ; relevons-le en un élément $g'$ de $\mathsf{Gal}(L/k)$. Par définition de $\mathsf{I}_{\mathsf {grad}}(F/k)$, l'action de $g'$ sur $\red{F}$ est triviale. Comme $\mathsf{Gal}(L/F)$ se surjecte sur $\mathsf{Gal} (\red{L}/\red{F})$, il existe $h$ appartenant à $\mathsf{Gal}(L/F)$ tel que $g'h^{-1}$ agisse trivialement sur $\red{L}$. Par construction, $g'h^{-1}$ est un élément de $\mathsf{I}_{\mathsf {grad}}(L/F)$ qui relève $g$.~$\Box$

\medskip
\deux{psylow} {\bf Proposition.} {\em Soit $L$ une extension finie galoisienne de $k$. Notons $\mathsf{I}(L/k)$ le noyau de $\mathsf{Gal}(L/k)\to \mathsf{Gal}(\red{L}_1/\red{k}_1)$. Le groupe $\mathsf{I}_{\mathsf {grad}}(L/k)$ est l'unique $p$-sous-groupe de Sylow de $\mathsf{I}(L/k)$.}

\medskip
{\em Démonstration.} Montrons tout d'abord que  $\mathsf{I}_{\mathsf {grad}}(L/k)$ est un $p$-groupe. Soit $g$ un élément de $\mathsf{I}_{\mathsf {grad}}(L/k)$ dont l'ordre $q$ est premier à $p$, et soit $E$ le sous-corps de $L$ formé des éléments invariants sous $g$. Soit $x$ un élément de $L$ ; l'élément $y=x-\mathsf{Tr}_{L/E}(x)/q$ de $L$ est de trace nulle sur $E$ ; on a donc $\sum g^{i}(y)=0$, soit encore $qy-\sum (g^{i}(y)-y)=0$. Comme $q$ est premier à $p$, on a $|qy|=|y|$ ; comme $g$ appartient à $\mathsf{I}_{\mathsf {grad}}(L/k)$, on a $|g^{i}(y)-y|<|y|$ pour tout $i$ si $y\neq 0$ ; par conséquent $y=0$ et $x\in E$. On vient de montrer que $E$ coïncide avec $L$, et partant que $g=\mathsf{Id}$ ; le groupe $\mathsf{I}_{\mathsf {grad}}(L/k)$ est donc bien un $p$-groupe.

\medskip
Il reste à prouver que $\mathsf{I}(L/k)/\mathsf{I}_{\mathsf {grad}}(L/k)$ est d'ordre premier à $p$ ; or en vertu de la proposition\ref{redquasigal}, ce quotient s'identifie à $\mathsf{Hom}(|L^{*}|/|k^{*}|,\red{L}_1^{*})$ ; on conclut en remarquant que la torsion de $\red{L}_1^{*}$ est première à $p$.~$\Box$ 

\medskip
\deux{maxmodram} {\bf Proposition.} {\em Soit $L$ une extension finie galoisienne de $k$ ; désignons par $F$ le sous-corps de $L$ formé des invariants sous $\mathsf{I}_{\mathsf {grad}}(L/k)$. On a l'égalité $[\red{F}:\red{k}]=[F:k]$. Par ailleurs, $\red{F}$ est une extension galoisienne de $\red{k}$, et la surjection naturelle $\mathsf{Gal}(\red{L}/\red{k})\to \mathsf{Gal}(\red{F}/\red{k})$ est bijective ; autrement dit, le corps $\red{F}$ est la fermeture séparable de $\red{k}$ dans $\red{L}$.}

\medskip
{\em Démonstration.} Comme $\mathsf{I}_{\mathsf {grad}}(L/k)$ est un sous-groupe distingué de $\mathsf G$ (en tant que noyau d'un homomorphisme), $F$ est une extension galoisienne de $k$. Comme $\mathsf{I}_{\mathsf {grad}}(L/k)$ agit trivialement sur $F$, la suite exacte du lemme \ref{wcompatible} assure que $\mathsf{I}_{\mathsf {grad}}(F/k)$ est trivial, et donc que $\mathsf{Gal}(F/k)\to\mathsf{Gal}(\red{F}/\red{k})$ est bijective. Des égalités et inégalités $$[F:k]=\sharp\mathsf{Gal}(F/k)=\sharp\mathsf{Gal}(\red{F}/\red{k})\leq [\red{F}:\red{k}]\leq [F:k],$$ il vient $[F:k]=\sharp\mathsf{Gal}(\red{F}/\red{k})=[\red{F}:\red{k}]$.

\medskip
Il en découle que $\red{F}$ est une extension galoisienne de $\red{k}$. D'autre part, $\mathsf{Gal}({\red{L}/\red{k}})$ s'identifie au quotient de $\mathsf {Gal}(L/k)$ par $\mathsf{I}_{\mathsf {grad}}(L/k)$, soit encore à $\mathsf{Gal}(F/k)$, lequel est naturellement isomorphe à $\mathsf{Gal}(\red{F}/\red{k})$ d'après ce qui précède. La surjection $\mathsf{Gal}(\red{L}/\red{k})\to\mathsf{Gal}(\red{F}/\red{k})$ est en conséquence bijective.~$\Box$ 

\medskip
\deux{extmodram} {\bf Proposition.} {\em Soit $L$ une extension finie galoisienne de $k$ et soit $E$ un sous-corps de $L$. Soit $F$ le sous-corps de $L$ formé des invariants sous $\mathsf{I}_{\mathsf {grad}}(L/k)$. Les propositions suivantes sont équivalentes :

\medskip
$i)$ $E$ est inclus dans $F$ ;

\medskip
$ii)$ $\red{E}$ est séparable sur $\red{k}$ et $[\red{E}:\red{k}]=[E:k]$.

\medskip
$iii)$ $[\red{E}:\red{k}]=[E:k]$, $\red{E}$ est séparable sur $\red{k}$ et $|E^{*}|/|k^{*}|$ est d'ordre premier à $p$.}

\medskip
{\em Démonstration.} L'équivalence de $ii)$ et $iii)$ résulte de la proposition~\ref{ef}. Montrons que $i)\Rightarrow ii)$. Supposons donc $E\subset F$. Dans ce cas, $\red{E}\subset \red{F}$, et $\red{E}$ est donc séparable sur $\red{k}$ par la proposition \ref{maxmodram}. Par ailleurs, comme $\mathsf{I}_{\mathsf {grad}}(L/k)$ fixe les éléments de $F$, on a l'égalité $\mathsf{I}_{\mathsf {grad}}(L/E)=\mathsf{I}_{\mathsf {grad}}(L/k)$; la proposition \ref{maxmodram} assure alors que $[\red{F}:\red{E}]=[F:E]$ ; comme  elle garantit également que $[\red{F}:\red{k}]=[F:k]$, on a bien $[\red{E}:\red{k}]=[E:k]$.

\medskip
Montrons maintenant que $ii)\Rightarrow i)$. Supposons donc que $E$ satisfait $ii)$ ; on a alors $\red E\subset \red F$ en vertu de la proposition \ref{maxmodram}. Soit $\xi$ un élément homogène non nul de $\gred{E}$, soit $r$ son degré et soit $R\in \gred{k}[X/r]$ son polynôme minimal sur $\red{k}$. Soit $\cal R$ un relèvement admissible de $R$ dans $k[X]$. Comme $R$ est séparable, le lemme de Hensel gradué assure que ${\cal R}$ possède une unique racine $x$ dans $E$ qui relève $\xi$, qu'il possède une unique racine $y$ dans $F$ qui relève $\xi$, et une unique racine $z$ dans $L$ qui relève $\xi$ ; en conséquence, $x=y=z$ et $\xi \in \red{(E\cap F)}$ ; il en découle que $\red{(E\cap F)}=\red{E}$.

\medskip
Dès lors $[E:k]=[\red{E}:\red{k}]=[\red{(E\cap F)}:k]=[(E\cap F):k],$ la première égalité étant satisfaite par hypothèse, la seconde parce que $\red{(E\cap F)}=\red{E}$, et la troisième en vertu de l'implication $i)\Rightarrow ii)$ déjà établie ; on en déduit que $E\subset F$.~$\Box$ 

\medskip
\deux{rempremp} {\bf Remarque.} Soit $E$ une sous-extension de $L$ telle que $[E:k]$ soit premier à $p$. Le groupe $\mathsf{Gal}(L/E)$ contient alors au moins un $p$-sous-groupe de Sylow $\mathsf S$ de $\mathsf{Gal}(L/k)$ ; comme $\mathsf{I}_{\mathsf {grad}}(L/k)$ est un $p$-sous-groupe distingué de $\mathsf{Gal}(L/k)$, il est contenu dans $\mathsf S$, et donc dans $\mathsf{Gal}(L/E)$ ; en conséquence, $E\subset F$.

\subsection*{Passage à la limite}

De la section précédente se déduisent aisément les faits qui suivent. 

\medskip
\deux{ksep} Soit $k^{s}$ une clôture séparable de $k$. On désigne par $\mathsf{I}_{\mathsf {grad}}(k^{s}/k)$ la limite projective des $\mathsf{I}_{\mathsf {grad}}(L/k)$, où $L$ parcourt la famille des sous-extensions finies et galoisiennes de $k^{s}$ ; on pose $M=(k^{s})^{\mathsf{I}_{\mathsf {grad}}(k^{s}/k)}$. Le corps $\red{M}$ est la fermeture séparable de $\red{k}$ dans $\red{k^{s}}$ et $\mathsf {Gal}(M/k)\simeq \mathsf{Gal}(\red{M}/\red{k})\simeq\mathsf{Gal}(\red{k^{s}}/\red{k})$ ; si l'on note $\Pi(k^{s}/k)$ ce dernier groupe, il s'insère dans une suite exacte naturelle $$1\to \mathsf{Hom}(|(k^{s})^{*}|/|k^{*}|,\red{(k^{s})}^{*})\to\Pi(k^{s}/k)\to\mathsf {Gal}(\red{k^{s}_1}/\red{k}_1)\to 1. $$

\medskip
On dit qu'une $k$-algèbre étale $E$ est {\em modérément ramifiée} si $E\otimes_{k}M$ est un produit fini de copies de $M$ ; l'algèbre $E$ est modérément ramifiée si et seulement si $\red{E}$ est une $\red{k}$-algèbre étale de dimension égale à $[E:k]$ (la structure d'algèbre normée sur $E$ a été définie au \ref{hypkhens}) ; la catégorie des $k$-algèbres étales modérément ramifiées est équivalente à celle des $\Pi(k^{s}/k)$-ensembles finis ; le foncteur $E\mapsto \red{E}$ établit une équivalence entre la catégorie des $k$-algèbres étales modérément ramifiées et celle des $\red{k}$-algèbres étales déployées par $\red{M}$.

\medskip
On peut améliorer cette dernière assertion. 

\medskip
\deux{deploye} {\bf Lemme.} {\em Toute $\red{k}$-algèbre étale est déployée par $\red{M}$.}

\medskip
{\em Démonstration.} Il suffit de montrer que $\red{M}$ est séparablement clos. Soit $r$ un élément de $\Gamma$ et soit $R$ un élément unitaire, homogène et séparable de $\red{M}[X/r]$ ; soit ${\cal R}$ un relèvement admissible de $R$ dans $M[X]$. Comme $R$ est séparable, $\cal R$ l'est également (l'image de son discriminant dans $\red{M}_{r}$ coïncide avec le discriminant de $R$) ; en conséquence, il a une racine $x$ dans $k^{s}$. Par construction, $\red{x}_{r}$ est une racine de $R$ dans $\red{k^{s}}$, qui est séparable sur $\red{k}$ et vit de ce fait dans $\red{M}$.~$\Box$ 

\medskip
\deux{conclu} Ainsi, {\em $E\mapsto \red{E}$ établit une équivalence entre la catégorie des $k$-algèbres étales modérément ramifiées et celle des $\red{k}$-algèbres étales}. Chacune de ces deux catégories est par ailleurs naturellement équivalente à celle des $\Pi(k^{s}/k)$-ensembles finis, le diagramme que l'on imagine étant (essentiellement) commutatif. 

\medskip
\deux{extcorps} Soit $\KK$ une extension $\Gamma$-valuée hensélienne de $k$ et soit $\KK^{s}$ une clôture séparable de $\KK$ ; soit $k^{s}$ la fermeture séparable de $k$ dans $\KK^{s}$. Il résulte des définitions que $\mathsf{Gal}(\KK^{s}/\KK)\to \mathsf{Gal}(k^{s}/k)$ envoie $ \mathsf{I}_{\mathsf {grad}}(\KK^{s}/\KK)$ dans $\mathsf{I}_{\mathsf {grad}}(k^{s}/k)$, et induit en particulier une flèche $\Pi(\KK^{s}/\KK)\to \Pi(k^{s}/k)$. Il en découle que si $E$ est une $k$-algèbre étale modérément ramifiée, alors $\KK\otimes_{k}E$ est une $\KK$-algèbre étale modérément ramifiée ; si ${\cal E}$ désigne le $\Pi(k^{s}/k)$-ensemble correspondant à $E$, alors le $\Pi(\KK^{s}/\KK)$-ensemble qui correspond à $\KK\otimes_{k}E$ est simplement celui déduit de $\cal E$ {\em via} la flèche $\Pi(\KK^{s}/\KK)\to \Pi(k^{s}/k)$. En utilisant le~\ref{conclu} ci-dessus, on voit que $\red{(\KK\otimes_{k}E)}$ s'identifie à $\red{\KK}\otimes_{\red{k}}\red{E}$. 

\section{Les polydisques virtuels}\label{PSEU}

\subsection*{La réduction d'un polydisque}

À partir de maintenant, le groupe ordonné $\Gamma$ est égal à $\RR^{*}_{+}$, et $k$ désigne un corps ultramétrique {\em complet}.

\medskip
\deux{redalg]} Soit $X$ un espace $k$-analytique au sens de Berkovich (\cite{brk1} et \cite{brk2} ; le choix d'une théorie géométrique précise est en réalité sans grande importance ici, {\em cf.} l'introduction) ; soit $\cal A$ l'algèbre des fonctions analytiques sur $X$. Pour tout réel $r$ strictement positif, on note ${\cal A}\lr$ (resp. ${\cal A}\str$) le sous-anneau de $\cal A$ formé des éléments $a$ tels que $||a||_{\infty}\leq r$ (resp. $||a||_{\infty}<r$) ; on désigne par ${\cal A}\hrr$ le complété de ${\cal A}\str$, pour la topologie de la convergence en norme sur tout domaine affinoïde, qui n'est autre que la topologie de la convergence uniforme sur tout compact dès que $X$ est réduit. Remarquons que si $X$ est compact, ${\cal A}\hrr={\cal A}\str$. On pose $\red{\cal A}=\bigoplus\limits_{r} {\cal A}\lr/{\cal A}\hrr$ ; c'est un anneau gradué. Si $a\in {\cal A}$ et si $r$ est un réel strictement positif supérieur ou égal à $||a||_{\infty}$, on notera $\red{a}_{r}$ l'image de $a$ dans $\red{\cal A}_{r}$. Si $r=||a||_{\infty}$, on écrira simplement $\red{a}$. La formation de $\red{\cal A}$ est clairement fonctorielle en $X$. 

\medskip
{\em Remarque.} Si $X$ est compact ({\em e.g}. affinoïde), alors $\cal A$ est une algèbre de Banach et $\red{\cal A}$ est sa réduction graduée au sens défini par Temkin.

\medskip
\deux{casdisque} {\bf Un exemple : le cas d'un polydisque.} Soit ${\bf r}=(r_{1},\ldots,r_{n})\in(\RR^{*}_{+})^{n}$, soit $X$ le $k$-polydisque {\em ouvert} de polyrayon $(r_{1},\ldots,r_{n})$ et soit $\cal A$ son algèbre des fonctions analytiques ; on désigne par $T_{1},\ldots,T_{n}$ les fonctions coordonnées sur $X$. 

\medskip
\trois{descarlarge} Soit $r$ un réel strictement positif. 

\medskip
$\bullet$ ${\cal A}\lr$ est l'ensemble des séries $\sum a_I {\bf T}^I$ telles que $|a_I|{\bf r}^I\leq r$ pour tout $I$ ; 

$\bullet$ ${\cal A}\str$ est l'ensemble des séries $\sum a_I {\bf T}^I$  pour lesquelles il existe $s<r$ tel que $|a_I|{\bf r}^I\leq s$ pour tout $I$ ; 

$\bullet$ ${\cal A}\hrr$ est l'ensemble des séries $\sum a_I {\bf T}^I$  telles que $|a_I|{\bf r}^I<r$ pour tout $I$. 

\medskip
\trois{descpolyresi} Quitte à réordonner les $r_{i}$ et les $T_{i}$, on peut supposer qu'il existe $j$ tel que $r_{i}$ soit de torsion modulo $|k^{*}|$ si et seulement si $i\leq j$. Pour tout $i\leq j$, notons $n_{i}$ l'ordre de $r_{i}$ modulo $|k^{*}|$ et choisissons $\lambda_{i}\in k^{*}$ tel que $r_{i}^{n_{i}}=|\lambda_{i}|$.

\medskip
Il découle du \ref{descarlarge} ci-dessus qu'il existe un isomorphisme de $\red{k}$-algèbres graduées entre $\red{\cal A}$ et $$\gred{k}[[\tau_{1},\ldots,\tau_{j}]]\; [r_{1}^{-1}t_{1},\ldots,r_{n}^{-1}t_{n}]/\;(\red{\lambda_{1}}^{-1}t_{1}^{n_{1}}-\tau_{1},\ldots,\red{\lambda_{j}}^{-1}t_{j}^{n_{j}}-\tau_{j}\;)$$ modulo lequel $\red{T_{i}}=t_{i}$ pour tout $i\in \{1,\ldots n\}$. 

\medskip
\trois{intrinspolyr} Le monoïde constitué par les degrés des éléments homogènes non nuls de $\red{\cal A}$ est égal à $|k^*|.r_1^{\NN}.r_2^{\NN}.\ldots r_n ^{\NN}$ ; ce dernier ne dépend donc que des propriétés intrinsèques ({\em i.e.} indépendantes de son plongement dans $\Aff^{n,{\rm an}}_k$) de l'espace $k$-analytique $X$. 

\medskip
\trois{anngradcom} L'anneau $\red{\cal A}$ est un anneau gradué intègre et {\em local}, dont l'idéal maximal $\got{m}$ est engendré par $(t_{1},\ldots,t_{n})$. Pour tout entier $d$, le $\red k$-espace vectoriel gradué $\red {\got m}^d/\red{\got m}^{d+1}$ (resp. $\red{\cal A}/\got m^{d+1}$) s'identifie au sous espace vectoriel gradué de $\red k [r_{1}^{-1}t_{1},\ldots,r_{n}^{-1}t_{n}]$ formé des polynômes qui sont homogènes de degré $d$ (resp. de degré $d$), {\em où le degré et l'homogénéité sont pris ici dans leur sens usuel, c'est-à-dire monomial.}  Donnons deux conséquences de ces faits : 

\medskip
$\bullet$ pour tout entier $d$, la flèche naturelle $${\rm S}^d(\got m/\got m^2)\to \got m^d/\got m^{d+1}$$ est un isomorphisme ; 

$\bullet$ l'anneau $\red{\cal A}$ est un anneau gradué local {\em complet}, ce qui signifie que pour tout $r\in\RR^{*}_{+}$, la flèche naturelle $\red{\cal A}_{r}\to \lim\limits_{\leftarrow} (\red{\cal A}/\got{m}^{d})_{r}$ est un isomorphisme. 

\medskip
\trois{resgraddiscextsc} Si $F$ est une extension complète quelconque de $k$, on déduit des descriptions explicites données ci-dessus (que l'on utilise sur $k$ et sur $F$) : que l'anneau des fonctions analytiques sur $X_F$ s'identifie naturellement à $F\hotimes_k{\cal A}$ (la complétion étant à prendre au sens de la convergence uniforme sur tout compact) ; et que $\red{(F\hotimes_{k}{\cal A})}$ s'identifie au complété de $\red{F}\otimes_{\red{k}}\red{\cal A}$ en son idéal homogène maximal $(t_{1},\ldots,t_{n})$. 

\medskip
\deux{systcoord} {\bf Définition.} {\em On conserve les notations $X,{\cal A}$, etc. introduites ci-dessus.} On dira qu'une famille $(f_{1},\ldots,f_{n})$ de fonctions sur $X$ est un {\em système de coordonnées} si les $f_{i}$ induisent un isomorphisme entre $X$ et un polydisque ouvert de $\Aff^{n,{\rm an}}_{k}$ centré en l'origine. Si c'est le cas et si l'on pose $s_i=||f_i||_{\infty}$ pour tout $i$ le polydisque ouvert en question est celui de polyrayon $(s_1,\ldots, s_n)$ ; remarquons qu'en vertu de \ref{intrinspolyr} on a alors $$|k^*|r_1^{\NN}r_2^{\NN}\ldots r_n ^{\NN}=|k^*|s_1^{\NN}s_2^{\NN}\ldots s_n ^{\NN}.$$

\medskip
\deux{lemmesystcoord} {\bf Proposition.} {\em Soit $(f_{1},\ldots,f_{n})$ une famille de fonctions bornées sur $X$. Les propositions suivantes sont équivalentes :

\medskip
\begin{itemize}

\item[$i)$] les $f_{i}$ forment un système de coordonnées de $X$ ; 
\item[$ii)$] les $\red{f_{i}}$ appartiennent à l'idéal maximal $\got{m}$ de $\red{\cal A}$ et en constituent une famille génératrice ; 
\item[$iii)$] les $\red{f_{i}}$ appartiennent à $\got{m}$ et leurs images dans $\got{m}/\got{m}^{2}$ en constituent une base (en tant que $\red{k}$-espace vectoriel gradué).
\end{itemize}}

\medskip
{\em Démonstration}. L'équivalence entre $ii)$ et $iii)$ résulte du lemme de Nakayama gradué ; l'implication $i)\Rightarrow ii)$ résulte du \ref{descarlarge} ; on suppose maintenant que $ii)$ est vraie, et l'on cherche à établir $i)$ ; pour tout $i$, on désigne par $s_{i}$ le rayon spectral de $f_{i}$. Soit $Y$ le polydisque ouvert de polyrayon $(s_1,\ldots,s_n)$ et soit $\phi$ le morphisme $X\to Y$ induit par les $f_i$ ; nous allons montrer que $\phi$ est un isomorphisme. 

\medskip
\trois{egalitcoeff} Écrivons $f_i=\sum a_{i,I}{\bf T}^I$ pour tout $i$ ; pour tout $j$ compris entre $1$ et $n$ on note $e_j$ le multi-indice $(\delta_{j\ell})_{\ell}$. Il découle de nos hypothèses que $|a_i,I|{\bf r}^I\leq \rho_i$ pour tout $(i,I)$, que $|a_{(i,0)}|<\rho_i$ pour tout $i$ et que le déterminant $\left|\red{a_{i,e_j}}_{\rho_ir_j^{-1}}\right|_{(i,j)}$ est un élément non nul de $\red k_{R}$, où $R=\left(\prod \rho_i\right)\left(\prod r_j\right)^{-1}$. 

\medskip
\trois{translatrho} La translation $\tau$ par $(-a_{1,0}, -a_{2,0},\ldots, -a_{n,0})$ induit en vertu de \ref{egalitcoeff} un automorphisme de $Y$ ; quitte à remplacer $\phi$ par $\tau\circ \phi$, on peut désormais supposer que $a_{i,0}=0$ pour tout $i$. 

\medskip
\trois{reversform} Il résulte de \ref{egalitcoeff} que le déterminant $\left|a_{i,e_j}\right|_{(i,j)}$ est de valeur absolue égale à $R$ ; il est en particulier non nul, ce qui signifie que $(f_1,\ldots,f_n)$ est une base de $\got n/\got n^2$, où $\got n$ est l'idéal maximal de $k[[{\bf T}]]$.  Par conséquent, le morphisme continu de $k[[{\bf T}]]$ dans lui-même induit par les $f_i$ est un isomorphisme ; cela entraîne en particulier l'existence pour tout $i$ d'une unique série $g_i\in \got n$ telle que $g_i(f_i)=T^i$. Nous allons montrer que chacune des $g_i$ converge sur $Y$, et que la fonction qu'elle y définit est de norme spectrale majorée par $r_i$. Le $n$-uplet $(g_1,\ldots, g_n)$ définira ainsi un morphisme $\psi: Y\to X$, qui vérifiera nécessairement $\psi \circ \phi={\rm Id}_X$ et $\phi\circ \psi={\rm Id}_Y$, par un calcul formel dans $k[[{\bf T}]]$. 

\medskip
\trois{bornegi} Écrivons $g_i=\sum b_{i,I}{\bf T}^I$ pour tout $I$. Fixons $i$ et donnons-nous un entier $d\geq 0$. Nous allons montrer que $|b_{i,I}|{\bf s}^I\leq r_i$ pour tout multi-indice $I$ tel que $|I|\leq d$, ce qui permettra de conclure. 

\medskip
Par définition de la famille $(g_1,\ldots,g_n)$ on a $$\sum_{|I|\leq d}b_{i,I}{\bf f}^I=T_i+h,$$ où la série $h$ ne comporte que des monômes dont le multi-degré est de taille strictement supérieure à $d$. Soit $\rho$ le maximum des $|b_{i,I}|{\bf s}^I$ pour $|I|\leq d$. La norme spectrale de $\sum\limits _{|I|\leq d}b_{i,I}{\bf f}^I$, vue comme fonction sur $X$, est alors majorée par $\rho$, et il en va donc de même de celle de chacun des monômes de la série $T_i+h$. On a dans $\red{\cal A}$ l'égalité $$\sum_{|I|\leq d} \red{b_{i,I}}_{\rho{\bf s}^{-I}}\red {\bf f}_{\bf s}^I=\red{T_{i,\rho}}+\red h _\rho.$$ Supposons que l'on ait $\rho>r_i$ ; on aurait alors $\red{T_{i,_\rho}}=0$, et $$\sum\limits_{|I|\leq d} \red{b_{i,I}}_{\rho {\bf s}^{-I}}\red {\bf f}_{\bf s}^I=\red h _\rho.$$ Comme $h$ ne comprend que des monômes dont le multi-degré est de taille strictement supérieure à $d$, l'élément $\red h_\rho $ de $\red{\cal A}$ appartient à $\got m^{d+1}$. On déduit alors de \ref{anngradcom} et du fait que les $f_i$ forment une base de $\got m/\got m^2$ que $\red{b_{i,I}}_{\rho {\bf s}^{-I}}=0$ pour tout $I$ tel que $|I|\leq d$, ce qui contredit la définition de $\rho$ comme maximum des $|b_{i,I}|{\bf s}^I$ pour $|I|\leq d$.~$\Box$

\medskip
\deux{pseudo} {\bf Théorème.} {\em Soit $X$ un espace $k$-analytique. Supposons qu'il existe une extension finie, séparable et modérément ramifiée $L$ de $k$ telle que $X_{L}$ soit isomorphe à un $L$-polydisque ouvert de polyrayon $(r_1,\ldots,r_n)$ ; l'espace $X$ lui-même est alors isomorphe à un $k$-polydisque ouvert dont le polyrayon $(s_1,\ldots, s_n)$ est tel que $$|L^*|r_1^{\NN}r_2^{\NN}\ldots r_n ^{\NN}=|L^*|s_1^{\NN}s_2^{\NN}\ldots s_n ^{\NN}.$$ }

\medskip
{\em Démonstration.} Choisissons un $n$-uplet $(l_{1},\ldots,l_{n})$ d'éléments de $L^{*}$ tel que $(\red{l_{1}},\ldots,\red{l_{n}})$ soit une base de $\red{L}$ sur $\red{k}$. Soit $\KK$ une extension ultramétrique complète quelconque de $k$. Comme $L$ est une extension modérément ramifiée de $k$, la famille $(\red{l_{1}},\ldots,\red{l_{n}})$ est d'après le~\ref{extcorps} encore une base de $\red{(L\otimes_{k}\KK)}$ sur $\red{k}$. Si $(\lambda_{1},\ldots,\lambda_{n})\in \KK^{n}$, le lemme \ref{libregredk} assure alors que $||\underbrace{\sum \lambda_{i}l_{i}}_{\in \KK\otimes_{k}L}||_{\infty}=\max |\lambda_{i}|.|l_{i}|$.

\medskip
Soit $\cal A$ (resp. $\cal B$) l'algèbre des fonctions analytiques sur $X$ (resp. $X_L$). Si $V$ est un domaine affinoïde de $X$ d'algèbre ${\cal A}_V$ alors $V_L$ est un domaine affinoïde de $X$ d'algèbre $L\otimes_k{\cal A}_V\simeq \bigoplus {\cal A}_V.l_i$ ; il s'ensuit que la flèche naturelle $$L\otimes_k{\cal A}=\bigoplus {\cal A}.l_i\to \cal B$$ est un isomorphisme. 

Par la définition de la norme spectrale d'une fonction analytique, et d'après ce qui précède, on a pour tout $n$-uplet $(a_{1},\ldots,a_{n})$ d'éléments de ${\cal A}$ l'égalité $$||\underbrace{\sum a_{i}l_{i}}_{\in \cal B}||_{\infty}=\max ||a_{i}||_{\infty}|l_{i}|.$$

On en déduit immédiatement que le morphisme canonique $\red{L}\otimes_{\red{k}}\red{\cal A}\to \red{\cal B}$ est un isomorphisme. 

\medskip
Par hypothèse, $X_{L}$ est isomorphe au $L$-polydisque ouvert de rayon $(r_1,\ldots,r_n)$. Si $m$ désigne sa dimension, $\red{\cal B}$ est un anneau gradué intègre local (et complet), dont toute famille génératrice minimale de l'idéal maximal homogène $\got{m}$ comprend $m$ éléments. Puisque $\red{\cal B}\simeq \red{L}\otimes_{\red{k}}\red{\cal A}$, la flèche $\red{\cal A}\to \red{\cal B}$ est injective et entière ; la version graduée du {\em going up} assure alors que $\red{\cal A}$ est local, d'idéal maximal $\got{n}=\got{m}\cap \red{\cal A}$. 

\medskip
Comme $L$ est modérément ramifiée sur $k$, la $\red{k}$-algèbre $\red{L}$ est étale ({\em cf.}~\ref{conclu}). L'anneau $\red{\cal B}/\got{n}$ étant dès lors une $\red {\cal A}/\got{n}$-algèbre finie étale et locale, c'est un corps, ce qui montre que $\got{n}$ engendre $\got{m}$. Il existe donc $m$ éléments bornés $f_{1},\ldots,f_{m}$ dans $\cal A$ tels que $(\red{f_{1}},\ldots,\red{f_{m}})$ engendre $\got{m}$. 

\medskip
Pour tout $i$ notons $s_{i}$ le rayon spectral de $f_{i}$ ; soit $Y$ le $k$-polydisque ouvert de polyrayon $(s_1,\ldots,s_n)$ et soit $\phi$ le morphisme $X\to Y$ induit par les $f_i$ ; il résulte de la proposition \ref{lemmesystcoord} que $\phi_L: X_L\to Y_L$ est un isomorphisme, et de \ref{intrinspolyr} que $|L^*|r_1^{\NN}r_2^{\NN}\ldots r_n ^{\NN}=|L^*|s_1^{\NN}s_2^{\NN}\ldots s_n ^{\NN}.$.

\medskip
Il reste à s'assurer que $\phi$ est un isomorphisme. On peut en fait montrer en toute généralité que le fait d'être un isomorphisme descend par extension quelconques des scalaires ({\em cf.} par exemple \cite{cnrtmk}, th. 9.2) ; mais lorsque les espaces en jeu sont sans bord, ce qui est le cas ici (il est clair que $Y$ est sans bord, et pour $X$ cela résulte du fait que $X_L$ est sans bord et que le morphisme fini $X_L\to X$ est sans bord, et par ailleurs surjectif), la preuve se simplifie, et nous allons la donner pour la commodité du lecteur. 

\medskip
Remarquons pour commencer que comme les fibres de $\phi_{L}$ sont des singletons, celles de $\phi$ le sont aussi ; ceci entraîne notamment la bijectivité de l'application continue sous-jacente à $\phi$.  Les applications $\phi_L$ et $Y_L\to Y$ sont topologiquement propres et en particulier fermées (la première parce qu'elle est finie, la seconde en tant que flèche induite par une extension du corps de base). Il en résulte, en vertu de la surjectivité de $X_L\to X$, que $\phi$ est fermée. 

\medskip
Soit $x\in X$ et soit $y$ son image sur $Y$. Comme $x$ est isolé dans la fibre $\phi^{-1}(y)=\{x\}$ et comme $X$ est sans bord (ce qui entraîne que $\phi$ est sans bord), le morphisme $\phi$ est fini en $x$. Il existe donc un voisinage affinoïde $U$ de $x$ et un voisinage affinoïde $V$ de $y$ tel que $\phi$ induise un morphisme fini $U\to V$ ; l'application $\phi$ étant fermée, et $x$ étant le seul antécédent de $y$, il existe un voisinage affinoïde $V'$ de $y$ dans $V$ tel que $U':=\phi^{-1}(V')$ soit contenu dans $U$ ; la flèche $U'\to V'$ est alors finie. 

\medskip
Soient $\cal C$ et $\cal D$ les algèbres affinoïdes respectivement associées à $U'$ et $V'$. La flèche $\phi_L$ étant un isomorphisme, $U'_L\to V'_L$ est un isomorphisme ; par conséquent, ${\cal D}_L\to {\cal C}_L$ est un isomorphisme, et comme ${\cal C}$ est une ${\cal D}$-algèbre de Banach finie, on a ${\cal C}_L={\cal C}\otimes_{\cal D}{\cal D}_L$. Il résulte alors de ce qui précède et de la fidèle platitude de la $\cal D$-algèbre ${\cal D}_L$ que $\cal D\to \cal C$ est un isomorphisme ; autrement dit, $U'\to V'$ est un isomorphisme. 

\medskip
Le morphisme $\phi$ apparaît ainsi comme un isomorphisme local sur sa source, qui est par ailleurs ensemblistement bijectif ; en conséquence, $\phi$ est un isomorphisme. ~$\Box$ 

\medskip
\deux{corollnonramray1} {\bf Corollaire.} {\em Soit $X$ un espace $k$-analytique et soit $n\geq 0$. Supposons qu'il existe une extension finie séparable et non ramifiée $L$ de $k$ tel que $X_L$ soit isomorphe au polydisque unité ouvert de dimension $n$ sur $L$. L'espace $X$ est alors lui-même isomorphe au polydisque unité ouvert de dimension $n$.} 

\medskip
{\em Démonstration.} Il s'agit simplement d'un cas particulier du théorème \ref{pseudo} ci-dessus : celui où $|L^*|=|k|^*$ et où les $r_i$ sont tous égaux à $1$. Les $s_i$ appartiennent alors à $|k^*|$, et le $k$-polydisque de polyrayon $(s_1,\ldots,s_n)$ est donc isomorphe au $k$-polydisque unité, ce qui termine la démonstration.~$\Box$ 

\medskip
\deux{commentnonram} {\em Commentaire.} Le corollaire ci-dessus peut se prouver directement, par une démonstration analogue à celle du théorème \ref{pseudo}, mais un peu plus simple, dans la mesure où elle évite le recours à l'algèbre graduée : il suffit de reprendre la section \ref{PSEU} depuis le début en remplaçant partout les réductions graduées par les réductions classiques (autrement dit, par leurs parties homogènes de degré 1) et en exigeant, dans la définition d'un système de coordonnées, que le but soit un polydisque {\em unité}. 

Modulo ces modification, tous nos arguments se retranscrivent {\em verbatim}, à l'exception de ceux donnés au \ref{bornegi} : si l'on veut éviter, pour ce point précis, le recours aux réductions graduées, il faut remarquer que $\rho\in |k^*|$ et diviser, si $\rho>1$, les deux termes de l'égalité $\sum\limits_{|I|\leq d}b_{i,I}{\bf f}^I=T_i+h$ par un scalaire de valeur absolue égale à $\rho$ ; c'est la réduction (classique) de l'égalité obtenue qui conduit alors à une contradiction. 

\medskip
\deux{contrex} {\bf Contre-exemples dans le cas sauvagement ramifié}. Nous allons tout d'abord donner un contre-exemple «générique», que nous allons ensuite décliner en deux cas particuliers. Supposons que $k$ est de caractéristique mixte $(0,p)$ et que $\gred k$ n'est pas parfait ; choisissons $a\in k^*$ tel que $\red a$ ne soit pas une puissance $p$-ième dans $\red k$, et posons $r=\sqrt [p] {|a|}$. Soit $L$ la $k$-algèbre étale $k[X]/X^p-a$. Si $F$ est un corps quotient de $L$ alors $\red a$ est une puissance $p$-ième dans $\red F$, ce qui montre que $[\red F:\red k]\geq p$, et donc que $[F:k]\geq p$. Il s'ensuit que $L$ est un corps de degré $p$ sur $k$, et que $\red L$ s'identifie à $\red k [X/r]/(X^p-\red a)$. 

\medskip
Soit $V$ l'ouvert de $\Aff^{1,\rm an}_k$ défini par la condition $|T^p-a|<|a|$ ; nous allons montrer que $V_L$ est isomorphe à un $L$-disque ouvert, mais que $V$ n'est pas lui-même isomorphe à un $k$-disque ouvert. 

\medskip
{\em L'espace $V_L$ est isomorphe à un $L$-disque ouvert.} Nous allons plus précisément vérifier que c'est un $L$-disque ouvert de $\Aff^{1,\rm an}_L$. Par construction, $a$ possède une racine $p$-ième $\xi$ dans $L$ ; notons que $|\xi|=r$. Soit $F$ un corps complet algébriquement clos quelconque contenant $L$ et soient $\xi=\xi_1,\xi_2,\ldots,\xi_p$ les racines $p$-ièmes de $a$ dans $F$. Si $i$ est un entier compris entre $2$ et $p$ alors  $\xi_i=\xi\mu$ pour une certaine racine $p$-ième de l'unité $\mu$ ; comme $|1-\mu|<1$, on a $|\xi_i|=r$ et $|\xi_i-\xi|<r$. Soit $z$ un élément de $F$. On déduit de ce qui précède que les équivalences suivantes : 

\medskip
$\bullet$ $|z-\xi|<r\iff (\forall i\; |z-\xi_i|<r)$ ; 

\medskip 
$\bullet$ $|z-\xi|\geq r\iff (\forall i\; |z-\xi_i|\geq r)$.

\medskip
Par conséquent, $$|z-\xi|<r\iff \prod |z-\xi_i|<r^p\iff |z^p-a|<|a|.$$ Ceci valant pour toute extension complète algébriquement close de $k$, l'ouvert $V_L$ de $\Aff^{1,\rm an}_L$ peut être défini par l'inégalité $|T-\xi|<r$, et est donc un $L$-disque ouvert ; remarquons que comme $r=|\xi|\in |L|^*$, l'ouvert $V_L$ est même isomorphe au disque {\em unité} ouvert. 

\medskip
{\em L'espace $V$ n'est pas isomorphe à un $k$-disque ouvert.} En effet, soit $x\in V$ et soit $\tau$ l'élément $T(x)$ de $\hres(x)$. On a $|\tau^p-a|<|a|$ ; il s'ensuit que $|\tau|=r$ et que $\red\tau ^p=\red a$ ; par conséquent, $\red {\hres(x)}\neq \red k$, et $\hres(x)$ n'est donc pas égal à $k$. Ainsi, $V$ n'a pas de $k$-point, et n'est dès lors pas isomorphe à un $k$-disque ouvert. 

\medskip
\trois{exsansram} {\em Première déclinaison de l'exemple ci-dessus.} On suppose que le corps résiduel {\em classique} $\red k_1$ n'est pas parfait, et l'on choisit $a$ dans $k\zero$ tel que $\red a $ ne soit pas une puissance $p$-ième dans $\red k_1$. L'ouvert $V$ de $\Aff^{1,\rm an}_k$ défini par l'inégalité $|T^p-a|<1$ devient isomorphe au disque unité ouvert sur $k(\sqrt [p]a)$, mais n'a pas de $k$-point, et n'est en particulier pas lui-même isomorphe à un $k$-disque ouvert. 

\medskip
\trois{exsansram} {\em Seconde déclinaison de l'exemple ci-dessus.} L'ouvert $V$ de $\Aff^{1,\rm an}_{\QQ_p}$ défini par l'inégalité $|T^p-p|<|p|$ devient isomorphe au disque unité ouvert sur $\QQ_p(\sqrt[p]p)$, mais n'a pas de $\QQ_p$-point et n'est en particulier pas lui-même isomorphe à un $k$-disque ouvert. 

\deux{analogcompl} {\em Une question.} Le théorème principal de cet article appelle la question naturelle suivante. A-t-on un résultat analogue dans le monde archimédien ? Autrement dit, toute forme réelle d'un polydisque unité ouvert complexe est-elle déjà un polydisque unité ouvert réel ? Pour donner un sens précis à cet énoncé, il n'est pas besoin de développer une théorie spécifique des espaces analytiques réels (esquissée par Berkovich au début de \cite{brk1}) ; on peut le traduire, en termes tout à fait classiques, de la façon suivante : {\em soit $n\in \NN$ et soit $\sigma$ une involution anti-holomorphe de $$\DD_n:=\{(z_1,\ldots,z_n)\in \CC^n,|z_i|<1\;\mbox{pour tout}\;i \}\;$$ existe-t-il un difféomorphisme holomorphe $\phi : \DD_n\to \DD_n$ tel que $\phi^{-1}\circ \sigma\circ \phi$ coïncide avec la conjugaison complexe ?} On pourrait bien sûr aussi poser la même question en remplaçant $\DD_n$ par la boule unité ouverte euclidienne de dimension $n$.

Lorsque $n=1$, on peut déduire facilement de la description explicite des automorphismes du disque que la réponse à cette question est positive. L'auteur remercie Antoine Chambert-Loir, Charles Favre et Henri Guenancia d'avoir (indépendamment, et presque simultanément) attiré son attention sur ce fait ; il ignore pour le moment ce qu'il en est en dimension quelconque.


\begin{thebibliography}{bb}

\bibitem{brk1} {\sc V. Berkovich}, {\em Spectral theory and analytic
    geometry over non-archimedean fields}, Mathematical Surveys and
    Monographs {\bf 33}, AMS, Providence, RI, 1990.
\bibitem{brk2} {\sc V. Berkovich}, {\em \'Etale cohomology for
    non-archimedean analytic spaces}, Inst. Hautes Etudes
    Sci. Publ. Math. {\bf 78} (1993) 5-161. 
 \bibitem{cnrtmk} {\sc B. Conrad, M. Temkin}, {\em Descent for non-archimedean analytic spaces}, preprint. 
 \bibitem{tbs}, {\sc Tobias Schmidt}, {\em Every unramified form of the non-Archimedean closed disc is trivial}, preprint. 
 \bibitem{tmk2} {\sc M. Temkin}, {\em On local properties of non-Archimedean analytic spaces. II.},  Israel J. Math.  {\bf 140}  (2004), 1-27.

\end{thebibliography}
\end{document}